\documentclass[3p]{elsarticle}
\makeatletter
\def\ps@pprintTitle{%
   \let\@oddhead\@empty
   \let\@evenhead\@empty
   \let\@oddfoot\@empty
   \let\@evenfoot\@oddfoot
}
\makeatother
\usepackage{graphicx, subfig, tikz} 
\usepackage{color}
\usepackage{multicol}
\usepackage{amsthm, amsmath, amssymb, mathtools, eucal, cancel, nicefrac} 
\usepackage{upgreek} 
\usepackage[scr=boondox]{mathalpha}

\usepackage[bbgreekl]{mathbbol}
\DeclareSymbolFont{AMSb}{U}{msb}{m}{n}
\DeclareSymbolFontAlphabet{\mathbb}{AMSb} 
\DeclareMathAlphabet{\bbol}{U}{bbold}{m}{n} 

\theoremstyle{definition}
\newtheorem{oss}{Remark}[section]

\newcommand{\m}[1]{\mathbf{#1}}

\newcommand{\RR}{\mathbb{R}}
\newcommand{\NN}{\mathbb{N}}
\newcommand{\norm}[1]{\lVert#1\rVert}
\DeclareMathOperator{\supp}{supp}
\newcommand{\defeq}{\coloneq}
\newcommand{\cR}{\mathcal{R}}
\newcommand{\cT}{\mathcal{T}}
\newcommand{\ceil}[1]{\left\lceil#1\right\rceil}
\renewcommand{\div}{\nabla\cdot}
\newcommand{\grad}{\nabla}
\newcommand{\eps}{\varepsilon}

\let\sl\relax
\DeclareMathOperator{\sl}{SL}
\DeclareMathOperator{\dl}{DL}
\DeclareMathOperator{\arccoth}{arccoth}

\DeclareMathOperator{\dsdt}{d\emph{s}\,d\emph{t}}
\DeclareMathOperator{\ds}{d\emph{s}}
\DeclareMathOperator{\dt}{d\emph{t}}
\DeclareMathOperator{\dsigma}{d\sigma_{\pmb{x}}}
\DeclareMathOperator{\dx}{d\pmb{x}}

\begin{document}
\begin{frontmatter}
\title{A Smoothly Varying Quadrature Approach for 3D IgA-BEM
Discretizations: Application to Stokes Flow Simulations}

\author[unifi]{Cesare Bracco}
\ead{cesare.bracco@unifi.it}
\author[TV]{Francesco Patrizi\corref{cor1}}
\ead{patrizi@mat.uniroma2.it}
\author[unifi]{Alessandra Sestini}
\ead{alessandra.sestini@unifi.it}
\address[unifi]{Department of Mathematics and Computer Science ``Ulisse Dini'',\\ University of Florence, Viale Giovanni Battista Morgagni 67/a, 50134 Florence, Italy}
\address[TV]{Department of Mathematics,\\ University of Rome Tor Vergata,  Via della Ricerca Scientifica, 1, 00133 Rome, Italy}
\cortext[cor1]{Corresponding Author}

\begin{abstract}
We introduce a novel quadrature strategy for Isogeometric Analysis (IgA) boundary element discretizations, specifically tailored to collocation methods. Thanks to the dimensionality reduction and the natural handling of unbounded domains, boundary integral formulations are particularly appealing in the IgA framework. However, they require the evaluation of boundary integrals whose kernels exhibit singular or nearly singular behavior. Even when the kernel is not singular, its numerical evaluation becomes challenging whenever the integration region lies close to a collocation point. These integrals of polar and nearly singular functions represent the main computational difficulty of IgA-BEM and motivate the development of efficient and accurate quadrature rules.
Unlike traditional methods that classify integrals as singular, nearly singular, or regular, our approach employs a desingularizing change of variables that smoothly adapts to the physical distance from singularities in the boundary integral kernels. The transformation intensifies near the polar point and progressively weakens when integrating over portions of the domain that are farther from it, ultimately leaving the integrand unchanged in the limit of a diametrically opposed region. This automatic calibration enhances accuracy and robustness by eliminating the traditional classification step, to which the approximation quality is often highly sensitive. Moreover, integration is performed directly over B-spline supports rather than over individual elements, reducing computational cost, particularly for higher-degree splines. The proposed method is validated through boundary element benchmarks for the three dimensional Stokes problem, where we achieve excellent convergence rates.
\end{abstract}

\end{frontmatter}
\section{Introduction}
Boundary Element Methods (BEMs) offer an attractive alternative to domain-based discretization techniques for the numerical solution of partial differential equations whenever Boundary Integral Equations (BIEs) can be formulated, as, e.g., in the cases of the Laplace, Helmholtz, linear elasticity, and Stokes problems, see \cite{Stein} for an in-depth introduction. Their appeal lies primarily in the dimensionality reduction achieved by reformulating the problem over the boundary and in their natural suitability for problems defined on unbounded domains.

Within the framework of Isogeometric Analysis (IgA), these advantages become even more significant, since IgA-BEMs exploit the exact boundary geometry provided by CAD models while avoiding volumetric meshing, which is often the most time-consuming and challenging step in standard domain-based IgA approaches. Both collocation and Galerkin formulations can be employed, with collocation being particularly popular in three-dimensional settings due to its simpler assembly process.

While dimensionality reduction offers clear benefits in the discretization phase, the boundary element approach also introduces two major challenges: the formation of dense linear systems and the treatment of singular and nearly singular integrals. The first issue can be efficiently mitigated through fast algorithms such as kernel-independent Fast Multipole Methods \cite{Ying04} or hierarchical matrix techniques \cite{Hackbush}, equipped with adaptive cross approximation \cite{Bebendorf}, as demonstrated, for instance, in \cite{sestini2}. In this work, we focus on the second challenge: the accurate numerical approximation of regular, (weakly) singular, and nearly singular integrals, particularly in the context of collocation IgA-BEMs.

Several approaches have been proposed to tackle this problem. In \cite{heltai}, analytical desingularization techniques were employed to regularize the kernel, enabling the use of standard Gaussian quadrature. Alternatively, the works \cite{sauertesi,sauer} combined Duffy transformations for singular cases, Gaussian quadrature for regular integrals, and modified Gauss--Legendre rules for nearly singular ones. The identification of nearly singular integrals in these methods has been simplified by considering only the integration over a ring of elements surrounding the one containing the collocation point.

We propose a new quadrature rule that unifies the treatment of regular and nearly singular integrals within a single adaptive framework, thereby avoiding the need for their preliminary classification. No analytical desingularization is required, a choice motivated by the kernel-dependent nature of such procedures and their additional complexity in multipatch configurations, see the discussion in \cite[Section~4.1]{sestini1}. Furthermore, we adopt a support-wise integration strategy, following the approach successfully applied to the Helmholtz equation in \cite{sestini1}, which significantly reduces the number of matrix accesses during assembly, compared to the element-wise strategies used in \cite{sauertesi,sauer,heltai}.

In the proposed quadrature rule, in the singular case, a subregion containing the singular point is isolated in the B-spline support and integration is done using a Duffy transformation in such subregion, ensuring accuracy. For the remaining regions, as well as for nonsingular supports, we employ an adaptive scheme that combines a smoothly varying change of variables, depending on the physical distance of the collocation point, generalizing the Telles transformation \cite{telles}, with a tensor polynomial interpolation of the kernel at Chebyshev points and direct spline integration formula \cite{manni}.
This combination yields a smoothly varying and adaptive quadrature scheme that automatically adjusts both number and distribution of quadrature nodes based on the distance between the collocation point and the integration support. As a result, the computational effort naturally decreases as the proximity to singularities diminishes, eliminating the need for explicit nearly singular classification while maintaining robustness and efficiency for all types of integrals.

The paper is organized as follows. Section \ref{sec:quadrule} constitutes the core of the work and presents the distance-calibrated change of variables that regularizes the integrand according to the proximity between the integration support and the singular point. This section also describes the numerical integration procedure for the non-singular case, which combines polynomial interpolation of the regularized kernels with direct support-wise B-spline integration. The treatment of singular integrals is also discussed: we show how to partition the B-spline support based on the location of the singular point and how to perform integration using our quadrature scheme, except in a small subregion containing the singularity, where a Duffy desingularization is employed.
In Section \ref{sec:stokes}, we recall the Stokes differential problem and outline the derivation of the associated boundary integral equation, which forms the basis for the BEM formulation. We therefore introduce the notation used in the numerical experiments and present two key integral identities involving the Stokes kernels. These identities are used primarily to check the quality of the proposed quadrature rule, as they allow us to isolate the integration error from the discretization error due to the approximation power of the discrete space.
Still in Section \ref{sec:stokes}, we detail the collocation IgA-BEM discretization process, establishing here as well some notation required for the numerical tests. Section \ref{sec:numericaltests} then presents the numerical results, which include tests of the proposed quadrature rule on the Stokes kernel identities as well as full collocation IgA-BEM simulations for standard Stokes benchmarks. The performance of the method is also compared with results available in the literature. Finally, Section \ref{sec:conclusion} draws the conclusions and outlines possible directions for future research.

\section{Smoothly varying quadrature formulas for weakly singular surface integrals}\label{sec:quadrule}
Focusing on the three dimensional setting and restricting attention to elliptic problems, we recall that for any boundary element discretization based on the collocation approach, the matrix assembly process requires the numerical approximation of integrals of the following type:
\begin{equation*} \label{refint}
\int_\Gamma K(\pmb{x}\,,\, \pmb{y}) B(\pmb{x})\, \dsigma
\end{equation*}
where $\Gamma \subseteq \RR^3$ is the finite boundary surface of the domain, $B\colon \Gamma \to \RR$ is a basis function of the adopted discretization space and $K\colon\Gamma \times \Gamma \to \RR,$ the kernel function associated with the underlying differential problem. Assuming $\Gamma$ is sufficiently smooth, i.e., endowed with a well-defined and continuous normal vector everywhere, the kernel $K$ exhibits a \textbf{weakly singular} behavior at $\pmb {x} = \pmb{y}$, namely
$$
K(\pmb{x}, \pmb{y}) \sim \frac{1}{\norm{\pmb{x} - \pmb{y}}}
\qquad \mbox{as  } \norm{\pmb{x} - \pmb{y}} \to 0$$
with $\pmb{x},\pmb{y} \in \Gamma$, see \cite{Stein} for details.
  We stress that the norm here considered is the Euclidean norm, representing the \emph{physical} distance between $\pmb{x}, \pmb{y}.$ As a consequence, even if numerical integration is performed on the surface $\Gamma,$ the singular behavior of $K$ depends on the distance in $\RR^3$ between $\pmb{x}$ and $\pmb{y}$, rather than on their geodesic distance on $\Gamma$. In particular, if $\Gamma$ is almost folded on itself, points that have a large geodesic distance could have a small physical distance. Therefore, when $\Gamma$ is parameterized as a multipatch surface, the points that are ``close" to $\pmb{y}$ do not necessarily belong to the same patch of $\pmb{y}$. 
  
  Usually, the basis function $B$ has local support and therefore the above integral is actually taking place just over $\supp(B)$. Furthermore, $B$ is the \emph{pushforward}, or \emph{lifting}, of a function defined over a parametric domain, which we can fix as $[0, 1]^2$ without loss of generality. More precisely, in the most general setting, where $\Gamma$ is a multipatch geometry composed of $M$ patches, the support of $B$ is contained in one of such patches. Assuming it is within the $m$--th patch, denoting the corresponding geometric mapping  as $\pmb{F}^{(m)}\colon[0, 1]^2 \to \RR^3$, it means that $B$ satisfies 
$$B(\pmb{x}) = B(\pmb{F}^{(m)}(\hat{\pmb{x}})) = \hat{B}(\hat{\pmb{x}})$$ with $\Gamma \ni \pmb{x} \eqqcolon \pmb{F}^{(m)}(\hat{\pmb{x}})$ for a $\hat{\pmb{x}} \in [0, 1]^2$ and $\hat{B}\colon [0, 1]^2 \to \RR$ the corresponding parametric basis function. In this context numerical quadrature is performed in the parametric domain through the \emph{pullback} operation, that is, the change of variables $\pmb{x} = \pmb{F}^{(m)}(\hat{\pmb{x}})$.
The integral we are analyzing therefore becomes
\begin{equation}\label{eq:BEMintegralpullback}
\int_{\supp(\hat{B})} K(\pmb{F}^{(m)}(\hat{\pmb{x}})\,,\, \pmb{y})\hat{B}(\hat{\pmb{x}})J^{(m)}(\hat{\pmb{x}})\, \text{d}\hat{\pmb{x}}
\end{equation}
with $J^{(m)}\colon[0, 1]^2 \to \RR$ the determinant of the Jacobian matrix of $\pmb{F}^{(m)}$.

We now distinguish two cases: \textbf{weakly singular integrals}, when $\pmb{y}$ belongs to $\supp(B)$ in $\Gamma$, and \textbf{non-singular integrals}, otherwise. 
The next Section \ref{sec:nonsing} addresses this latter class of integrals and represents the main contribution of the paper. Our approach is based on a tensor product bicubic change of variables that regularizes the integrand according to the physical distance between the integration support of $B$ and the singular point. The regularized kernel is then approximated through polynomial interpolation at Chebyshev nodes. The intensity of the regularizing transformation and the number of interpolation nodes are both selected automatically, according to such distance. It is important to emphasize that interpolation is not performed on $\hat{B}$ but only on the kernel and Jacobian factors. In case of B-spline discretizations, in the IgA-BEM context, this avoids potential issues arising from the limited inter-element regularity of $B$, which is only modified by the change of variables, preserving its tensor product spline structure. Thereby, we also enable direct spline integration and, consequently, support-wise integration. This allows for a function-by-function assembly process, rather than the conventional element-by-element approach, reducing memory access since each of such matrix entries is thereby visited only once.
Section \ref{sec:sing} presents the treatment of weakly singular integrals. In contrast to the non-singular case, the integration domain is first partitioned into at most five rectangular subregions, with the singular point $\pmb{y}$ lying in one of them. In this subregion, a Duffy transformation is applied to handle the singularity, while in all remaining subregions the same adaptive quadrature rule introduced for non-singular integrals is employed.
Finally, Section \ref{sec:NS} describes the distance dependent distribution and density of quadrature nodes. Although the average number of nodes per support is kept essentially constant, we recommend clustering and locally increasing the node density toward the singular point. This strategy, combined with the regularizing change of variables, promotes uniform accuracy across all supports.

\subsection{Non-singular integrals}\label{sec:nonsing}
When the polar point $\pmb{y}$  does not belong to $\supp(B),$ integration is not singular. However, the closer $\pmb{y}$ is to $\supp(B)$, the more challenging its numerical computation becomes. Thus, traditionally, non-singular integrals are classified into \emph{nearly singular} and \emph{regular} integrals. The distinction is based on either a distance threshold from the polar point or on the number of \emph{element rings} around the element containing $\pmb{y}$, where the first ring is made of the elements sharing an edge or a vertex with such element, the second ring by the elements sharing an edge or a vertex with the first ring, and so on and so forth. In both cases, an arbitrary setting of a parameter for the classification is required. Moreover, with the second criterion based on element rings, the separation of the integrals changes as we refine the discretization, with the region formed by the union of the elements in the integrals of the nearly singular class shrinking as we refine the underlying mesh. Furthermore, as we mentioned, elements that are physically ``close'' to $\pmb{y}$ but have a ``large'' geodesic distance from it would not be included in the nearly singular group when this second criterion is used. Once the classification is done, the standard approach is to use a special, more complicated and expensive, quadrature rule for the nearly singular integrals and a simpler and cheaper rule for the regular integrals. Therefore, the procedure is very sensitive to the integral classification, both in terms of accuracy and computational cost. 

With the purpose of avoiding such a classification, we propose a unique quadrature rule which takes as input, and continuously depends on, the physical distance $\delta$ between $\supp(B)$ and $\pmb{y}$, in order to ensure a correct consideration of the polar point position with respect to the domain of integration. This measure shall \emph{calibrate} how relevant $\pmb{y}$ is for the numerical integration: the further it is, the less we care about its position with respect to $\supp(B)$. In the limit, when $\delta$ is essentially equal to the diameter of $\Gamma$, $\text{diam}(\Gamma)$, we shall not care about it at all. Thereby, we make a \textbf{smoothly varying quadrature rule} and work around the classification of non-singular integrals in regular and nearly singular. The rule is composed of two ingredients: a change of variables and  polynomial interpolation. The change of variables shall regularize the weakly singular function $K$ to some extent according to the distance $\delta$: the smaller $\delta$, the stronger the transformation, while, as it becomes larger and approaches $\text{diam}(\Gamma)$, the transformation tends to the identity. We named this change of variables the \textbf{Distance Calibrated Telles (DCT) transformation}.

\subsubsection{Distance Calibrated Telles transformation}
\begin{figure}
    \centering
    \subfloat[]{\includegraphics[width=0.3\textwidth, page = 1]{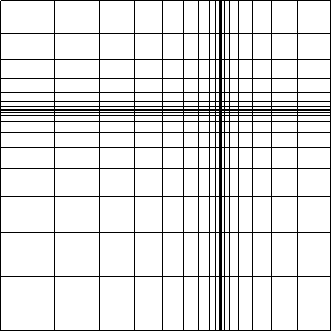}}\quad
    \subfloat[]{\includegraphics[width=0.3\textwidth, page = 2]{dct.pdf}}\quad
    \subfloat[]{\includegraphics[width=0.3\textwidth, page = 3]{dct.pdf}}\\
    \subfloat[]{\includegraphics[width=0.3\textwidth, page = 4]{dct.pdf}}\quad
    \subfloat[]{\includegraphics[width=0.3\textwidth, page = 5]{dct.pdf}}\quad
    \subfloat[]{\includegraphics[width=0.3\textwidth, page = 6]{dct.pdf}}
    \caption{Effect of DCT transformations on a uniform grid of isoparametric lines in $[-1, 1]^2$ for different values of the parameter $\alpha \in [0, 1]$ and a fixed point $\pmb{s} = \pmb{s}_\alpha$ in the domain, used for the definition of all the mappings. In (a) $\alpha = 0$ and we can see the effect of the tensor-product Telles transformation $\pmb{q}_0$. The isolines concentrate toward $\pmb{s}_\alpha$. In (b)--(e) we see the smooth decay in the intensity of the transformation, that is, on the deformation of the grid, for increasing values of $\alpha$, ultimately having no effect in (f) for $\alpha = 1$, corresponding to the identity map.}
    \label{fig:exDCTiso}
\end{figure}

The DCT transformation that we are going to introduce is a univariate change of variable for $\hat{x} \in [-1, 1]$. Therefore, we shall apply it in a tensor fashion when $\hat{\pmb{x}} \in [-1, 1]^2$ after a scaling of $\supp(\hat{B})$ to $[-1, 1]^2$. Let $\alpha = \alpha(\delta)$ be an increasing function of the physical distance $\delta$ with values in the range $[0, 1]$. We look for a family of polynomial functions $\{q_\alpha\}_{\alpha \in [0, 1]}$ such that, for every value of $\alpha$, $q_\alpha$ is invertible and preserves the integration interval, that is, $q_\alpha(\pm 1) = \pm 1$. Furthermore, in the limit cases, when $\alpha = 1$ we want $q_1$ to be the identity map and when $\alpha = 0$ we want $q_0$ to be a change of variables that desingularizes the integral, which has a polar singularity at $\hat{y} \in [-1, 1]$. For all the other values $\alpha \in (0, 1)$ we shall have a gradual transition between these two opposite situations, therefore generating a regularization spectrum provided by such changes of variables. A well established desingularization technique based on an invertible polynomial change of variables, preserving also the interval $[-1, 1]$, is the \emph{Telles transformation} \cite{telles}. The desingularization and invertibility requirements impose the existence of a point $s_0 \in [-1, 1]$ such that $q_0(s_0) = \hat{y}$, $q_0'(s_0) = 0$ and $q''(s_0) = 0$. In particular, when substituting the variables in the integral, the presence of the (univariate) Jacobian determinant, $q_0'$, together with the polynomial expression of $q_0$, ensures the well-posedness of the new integrand everywhere in $[-1, 1]$, even in $s_0$. The least degree polynomial that verifies the constraints of invertibility and desingularization has the following form:
\begin{equation}\label{eq:Telles}
q_0(s) = \hat y + \beta_0(s - s_0)^3
\end{equation}
for $\beta_0 \in \RR\setminus\{0\}$. Then by enforcing interval preservation we determine $\beta_0$ as
$$
\beta_0 = \frac{1}{1 + 3s_0^2}
$$
and $s_0 \in [-1, 1]$ as the solution of the following cubic equation:
$$
\hat{y}(1+3 s_0^2) = s_0 (3 + s_0^2)
$$
which can be solved exactly using Tartaglia's formula for determining the unique real root $s_0 \in [-1,1].$ \cite{cardano}. 
The DCT transformation generalizes the Telles change of variables by allowing a linear term in the expression:
\begin{equation}\label{eq:DCT}
q_\alpha(s) = \hat{y} + \alpha(s - s_\alpha) + \beta_\alpha(s - s_\alpha)^3,
\end{equation}
for $\alpha \in [0, 1]$. Again, the values of $\beta_\alpha$ and $s_\alpha$ are determined by imposing the interval preservation conditions $q_\alpha(\pm 1) = \pm 1$, which lead to
$$
\beta_\alpha = \frac{1 - \alpha}{1 + 3s_\alpha^2}
$$
and the following cubic equation in $s_\alpha$
$$
\hat{y}(1+3 s_\alpha^2) = s_\alpha [3-2\alpha + (2\alpha + 1)s_\alpha^2].
$$
which again can be solved using Tartaglia's formula \cite{cardano} or the more numerically stable formula in \cite{vignes} (needed especially as $\alpha$ approaches $1$). 
Note that the invertibility of $q_\alpha$ is guaranteed by the assumption that $\alpha$ ranges between $0$ and $1$, which implies also the positivity of $\beta_\alpha$, so that $q_\alpha$ is increasing monotone. Instead, we stress that $q_\alpha'(s_\alpha) = \alpha \neq 0$ for $\alpha > 0$. However, desingularization is not needed for $\alpha > 0$ because there will not be a polar point in the integration domain, as the distance $\delta$ is strictly positive. For $\alpha > 0$, that is, in case of a non-singular integral, and for a fixed direction, the special point $\hat{y}$ used for the construction of the DCT transformation shall be the coordinate in such direction of the point $\hat{\pmb{y}} \in [-1, 1]^2$ whose pushforward is the closest point in $\supp(B)$ to $\pmb{y}$.
Finally, we observe that the DCT transformation \eqref{eq:DCT} is consistent with the Telles transformation, that is, $q_0$ is indeed \eqref{eq:Telles}, and on the contrary $q_1$ is the identity mapping. In fact, for $\alpha = 1$, we have $\beta_1 = 0$ and the equation of $s_1$ becomes
$$
\hat{y}(1 + 3s_1^2) = s_1(1 + 3s_1^2)
$$
which leads to $s_1 = \hat{y}$. Therefore,
$$
q_1(s) = \hat{y} + (s - \hat{y}) = s
$$
i.e., the identity mapping. In Figure \ref{fig:exDCTiso} we show the effect of tensor product DCT transformations $\pmb{q}_\alpha\colon[-1, 1]^2 \to [-1, 1]^2$ in the bivariate setting on isoparametric lines in $[-1, 1]^2$, for the same function $\alpha$ in the two directions for simplicity. The extremal values $\alpha = 0$ (Telles transformation), $\alpha = 1$ (identity) and some intermediate DCT change of variables are shown, illustrating the smooth changing of intensity in the outcomes.

\begin{figure}
    \centering
    \includegraphics[width=0.35\linewidth]{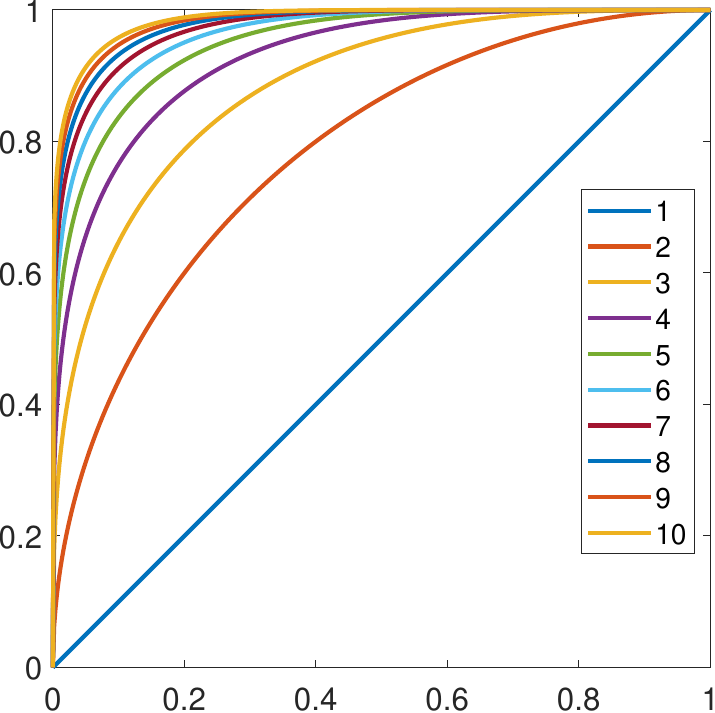}
    \caption{Functions belonging to the class \eqref{eq:alphaclass}, for  $\gamma = 1, \ldots, 10$. On the $x$ axis the normalized distance $\nicefrac{\delta}{\text{diam}(\Gamma)}$. The higher $\gamma$, the faster the function grows for small values of $\delta$ and, symmetrically, the slower it grows for large values of $\delta$.}
    \label{fig:alphaplot}
\end{figure}

As we mentioned, we shall have $\alpha$ as an increasing function of the distance $\delta$ and with range $[0, 1]$. Here, there is no preferred choice. For symmetry and beauty, we suggest picking $\alpha$ from this class of functions
\begin{equation}\label{eq:alphaclass}
\{\alpha_\gamma(\delta) \defeq (1 - (1 - \nicefrac{\delta}{\text{diam}(\Gamma)})^\gamma)^{\nicefrac{1}{\gamma}} \,:\, \gamma \in \NN \setminus\{0\}\},
\end{equation}
which also proved to work well in our tests. 
Some functions of this class are plotted in Figure \ref{fig:alphaplot} which shows that the higher $\gamma$, the faster $\alpha$ grows with $\delta$ and hence the sooner the relevance of the polar point for the numerical integration is lightened, as the regularization is mitigated.  

\subsubsection{Polynomial Interpolation \& Spline Direct Integration}
Once the integrand has been regularized with the DCT transformation, according to the distance of $\supp(B)$ from the polar point $\pmb{y}$,  an interpolatory numerical quadrature is performed to approximate the value of the integral. We can now employ the same rule for this purpose, regardless of the original difficulty of such integral, as the change of variables uniformizes the complexity of the problem.
The quadrature rule we are going to present is specific for spline discretizations in the so called IgA-BEM context, as we shall exploit B-spline properties. Hence, we now assume $\hat{B}$ to be a B-spline of bi-degree ${\bf d} = (d_1, d_2)$. A short introduction to B-splines and their features can be found in \cite{manni}. Let $\pmb{A}\colon[-1, 1]^2 \to \supp(\hat{B})$ be the affine change of variables we need for converting the domain of integration from $\supp(\hat{B})$ to $[-1, 1]^2$ in order to apply the DCT transformation. Then, if $\bbol{J}_{\pmb{A}}$ is the Jacobian matrix of $\pmb{A}$, it holds $\det \bbol{J}_{\pmb{A}} = \nicefrac{|\supp(\hat{B})|}{4}$. Defined now the composition of the changes of variables $\pmb{G}_\alpha \defeq \pmb{q}_\alpha \circ \pmb{A}$,  $\pmb{G}_\alpha\colon [-1, 1]^2 \to \supp(\hat{B})$, with $\pmb{G}_\alpha(s,t) = \hat{\pmb{x}}$, with $\pmb{q}_\alpha$ the tensor-product DCT transformation, we have that its Jacobian determinant is $\det \bbol{J}_{\pmb{G}_\alpha} \defeq \nicefrac{|\supp(\hat{B})|}{4}\det \bbol{J}_{\pmb{q}_\alpha}$ and the integral \eqref{eq:BEMintegralpullback} becomes 
$$
  \int_{[-1, 1]^2} K(\pmb{G}_\alpha(s, t)\,,\, \pmb{y})\, \hat{B}(\pmb{G}_\alpha(s, t))\, J^{(m)}(\pmb{G}_\alpha(s, t))\,\det\bbol{J}_{\pmb{G}_\alpha}(s, t)\, \dsdt.
$$ 
The change of variables $\pmb{G}_\alpha$ has univariate components, i.e., $\pmb{G}_\alpha(s, t) = (G_\alpha^s(s), G_\alpha^t(t))$ for $(s, t) \in [-1, 1]^2$, with $G_\alpha^s, G_\alpha^t$ polynomials of degree $3$.
The bivariate B-spline $\hat{B}$ is a separable function, that is, $\hat{B}(\hat{x}_1, \hat{x}_2) = \hat{B}^1(\hat{x}_1) \cdot \hat{B}^2(\hat{x}_2)$ for $\hat{x}_1, \hat{x}_2 \in [0, 1]$ with $\hat{B}^1, \hat{B}^2$ univariate B-splines of degrees $d_1, d_2$, respectively. Therefore, we have that $\hat{B} \circ \pmb{G}_\alpha = (\hat{B}^1 \circ G_\alpha^s) \cdot (\hat{B}^2 \circ G_\alpha^t)$ is a bivariate spline of bi-degree $(3d_1, 3d_2)$ defined in $[-1, 1]^2$ and with same regularity of $\hat{B}$. In order to lighten the notation, setting 
\begin{equation} \label{eq:Ktilde}
\tilde{K}_{\pmb{y}}(s, t) \defeq K(\pmb{G}_\alpha(s, t)\,,\,\pmb{y})\, J^{(m)}(\pmb{G}_\alpha(s, t))\, \det\bbol{J}_{\pmb{G}_\alpha}(s, t)
\end{equation}
for all $(s, t) \in [-1, 1]^2,$ we get the following compact expression for the integral \eqref{eq:BEMintegralpullback},
\begin{equation*}\label{final}
\int_{[-1, 1]^2} \tilde{K}_{\pmb{y}}(s, t) \,\hat{B}(\pmb{G}_\alpha(s, t))\,\dsdt.
\end{equation*}
We emphasize that the function $\tilde{K}_{\pmb{y}}$ in \eqref{eq:Ktilde} can be numerically integrated more easily than its original expression. This is because the kernel $K$ has been regularized through the reparametrization $\pmb{G}_\alpha$, which stabilizes its evaluation even when $K$ was originally defined near a singularity. Moreover, the Jacobian determinant of the patch map, $J^{(m)}$, is assumed to be smooth over the entire supports of the pushforward B-splines. In particular, when the mapping $F^{(m)}$ is a NURBS with internal knots, smoothness of $J^{(m)}$ is guaranteed in $\supp(B)$ by embedding the NURBS geometry knots into the discretization mesh with maximal multiplicity.

Given a pair of integers $(c_1, c_2)$, we approximate $\tilde{K}_{\pmb{y}}$ with the tensor product polynomial interpolant of bi-degree $(c_1, c_2)$ in the Chebyshev nodes $(s_i)_{i=0}^{c_1}$ and $(t_j)_{j=0}^{c_2}$, respectively. Let $\{\ell_i^s\}_{i=0}^{c_1}$ and $\{\ell_j^t\}_{j=0}^{c_2}$ be the (univariate) Lagrange polynomial bases with respect to such Chebyshev nodes in the two directions. Then 
\begin{equation*}
\begin{split}
\int_{[-1, 1]^2} \tilde{K}_{\pmb{y}}(s, t)\, \hat{B}(\pmb{G}_\alpha(s, t))\, \dsdt &\approx \int_{[-1, 1]^2} \sum_{j=0}^{c_2 }\sum_{i = 0}^{c_1} \tilde{K}_{\pmb{y}}(s_i, t_j) \ell_i^s(s)\ell_j^t(t)\, \hat{B}(\pmb{G}_\alpha(s, t))\, \dsdt\\
&= \int_{-1}^{1}\sum_{j=0}^{c_2} \left(\int_{-1}^1 \sum_{i=0}^{c_1 } \tilde{K}_{\pmb{y}}(s_i, t_j) \ell_i^s(s) \hat{B}^1(G_\alpha^s(s))\, \ds\right)\ell_j^t(t)\hat{B}^2(G_\alpha^t(t))\, \dt\\
&= \int_{-1}^{1}\sum_{j=0}^{c_2 } \left(\int_{-1}^1 L_j^s(s) \hat{B}^1(G_\alpha^s(s))\, \ds\right)\ell_j^t(t)\hat{B}^2(G_\alpha^t(t))\, \dt\\
&= \int_{-1}^1 \sum_{j=0}^{c_2} I_j\ell_j^t(t)\hat{B}^2(G_\alpha^t(t))\, \dt\\
& = \int_{-1}^1 L^t(t)\, \hat{B}^2(G_\alpha^t(t))\, \dt,
\end{split}
\end{equation*}
where we have set for all $j = 0, \ldots, c_2$
$$
L_j^s(s) \defeq \sum_{i=0}^{c_1} \tilde{K}_{\pmb{y}}(s_i, t_j) \ell_i^s(s), \quad I_j \defeq \int_{-1}^1 L_j^s(s)\, \hat{B}^1(G_\alpha^s(s))\, \ds, \quad\text{and}\quad
L^t(t) \defeq \sum_{j=0}^{c_2}I_j\ell_j^t(t).
$$
Therefore, we compute first $I_j$ for all $j = 0, \ldots, c_2$ and then the final integral of $L^t\cdot (B^2 \circ G_\alpha^t)$. Clearly, a mirrored procedure is possible when swapping the ordering of summation between $i$ and $j$ and the most computationally convenient option can be chosen. We focus first on the integrals $I_j$. The product of a polynomial $L_j^s$ of degree $c_1$ with a spline $(\hat{B}^1\circ G_\alpha^s)$ in $[-1, 1]$ is a univariate spline, that we shall call \textbf{product spline}, of degree $3d_1 + c_1$, whose breakpoints in $[-1, 1]$ are the anti-images via $(G_\alpha^s)^{-1}$ of the breakpoints of $\hat{B}^1$ with multiplicity raised by $2d_1 + c_1$ with respect to their multiplicity in the local knot vector of $\hat{B}^1$, in order to achieve the given smoothness. Namely, if $\mu_k$ is the multiplicity of the $k$--th breakpoint of $\hat{B}^1$, then its multiplicity in the spline product knot vector will be $\tilde{\mu}_k = \mu_k + 2d_1 + c_1$. Having established the knot vector where the product spline of degree $3d_1 + c_1$ is defined, we can consider the associated B-spline basis  $\{\tilde{B}_k\}_{k=1}^{\nu_1}$. We therefore have
$$
L_j^s(s)\hat{B}^1(G_\alpha^s(s)) = \sum_{k = 1}^{\nu_1} \lambda_k^{j, s} \tilde{B}_k(s).
$$
In order to find the coefficients $\pmb{\uplambda}^{j, s} = (\lambda_k^{j, s})_{k=1}^{\nu_1}$ we solve an interpolation (collocation) problem $\bbol{C}^s\pmb{\uplambda}^{j, s} = \m{v}^{j, s}$ with $\bbol{C}^s \defeq (\tilde{B}_{k_2}(\tilde{s}_{k_1}))_{k_1, k_2 = 1}^{\nu_1}$ where $\{\tilde{s}_k\}_{k=1}^{\nu_1}$ are the Greville points of the B-spline basis $\{\tilde{B}_k\}_{k=1}^{\nu_1}$ and $\text{v}_{k_1}^{j, s} \defeq L_j^s(\tilde{s}_{k_1})\hat{B}^1(G_\alpha^s(\tilde{s}_{k_1}))$. The integral of the B-spline $\tilde{B}_k$ can finally be computed exactly and directly through the well-known integration formula for B-splines, see e.g. \cite[Theorem 13]{manni}, i.e., if $\supp(\tilde{B}_k) = [\tilde{\xi}_k, \tilde{\xi}_{k + 3d_1 + c_1 + 1}]$ we have
$$
\int_{-1}^1 \tilde{B}_k(s)\ds = \int_{\supp(\tilde{B}_k)} \tilde{B}_k(s)\ds = \nicefrac{(\tilde{\xi}_{k+3d_1 + c_1 + 1} - \tilde{\xi}_k)}{(3d_1 + c_1 + 1)} \eqcolon c_k^s.
$$
Therefore, 
$$
I_j = \int_{-1}^1 L_j^s(s) B^1(G^s(s))\ds = \sum_{k=1}^{\nu_1} \lambda_k^{j, s} c_k^s.
$$
We then apply the same procedure to the other direction to compute the final integral exactly:
\begin{equation*}
\int_{-1}^1 L^t(t)\hat{B}^2(G_\alpha^t(t))\dt =  \int_{-1}^1 \sum_{\ell = 1}^{\nu_2} \lambda_\ell^t \tilde{B}_\ell(t)\dt = \sum_{\ell = 1}^{\nu_2} \lambda_\ell^t c_\ell^t. 
\end{equation*}
Again the coefficients $\pmb{\uplambda}^t = (\lambda_\ell^t)_{\ell = 1}^{\nu_2}$ of the representation of the product spline $L^t \cdot (\hat{B}^2\circ G_\alpha^t)$ in terms of the B-splines $\{\tilde{B}_\ell\}_{\ell = 1}^{\nu_2}$ of degree $3d_2 + c_2$ can be computed solving the interpolation problem $\bbol{C}^t\pmb{\uplambda}^t = \m{v}^t$ at the Greville points of such B-spline basis.

The proposed quadrature formula for non-singular integrals is of weighted interpolatory type. Specifically, the integrand in \eqref{eq:BEMintegralpullback} is the product of two functions: a regularized kernel, and the composition of a tensor product B-spline with the tensor- product polynomial reparametrization $\pmb{G}_{\alpha}$. The latter is left unchanged by the quadrature rule, i.e., it is not approximated by a polynomial, and therefore acts as a weight function.
Conversely, the regularized kernel is interpolated at suitably gridded points by a tensor product polynomial. The weighted nature of the quadrature rule can be mantained, which is particularly advantageous for integration over B-spline supports, while adopting a different type of approximation for the regularized kernel, for instance, tensor product uniform spline quasi-interpolation, as proposed in \cite{sestini1}, where $h$-refinement was mainly considered to improve the rule accuracy.
In the present work, however, we prefer polynomial interpolation at Chebyshev points, since $p$-refinement is generally more effective for integrals with smooth integrands. We recall that the kernel remains well-behaved even near the singularity, thanks to the preliminary DCT transformation, and thus it remains smooth in this situation as well, allowing the use of $p$-refinement with polynomial interpolation.
\subsection{Singular integrals}\label{sec:sing}
We now present a quadrature rule for integrating \eqref{eq:BEMintegralpullback} when $\pmb{y} \in \supp(B)$. 
Instead of applying also in this case Telles transformation, we propose a desingularization leading to a more regular function, therefore requiring a lower computational effort for integration. Furthermore, we observe that only less than $4(d_1 + 1)(d_2 + 1)$ integrals are singular because of the local linear independence of B-splines on the elements containing $\pmb{y}$ (which in the worst case scenario is an internal vertex of the mesh, shared by 4 elements). Hence, this quadrature rule is applied only to a fixed small number of the integrals and it only happens in a small portion of $\Gamma$. As we shall see, the procedure for these integrals requires to subdivide $\supp(\hat{B})$ in at most 5 regions, regardless of the degree $(d_1, d_2)$.  Only in one of these 5 regions a different approach is used with respect to the one adopted for the non-singular integrals. The procedure works as follows.

\begin{figure}[h!]
\centering
\subfloat[]{
\begin{tikzpicture}[scale = 1]
\draw (3, 2) rectangle (4, 4);
\draw (1, 3) rectangle (3, 4);
\draw (1, 1) rectangle (2, 3);
\draw (2, 1) rectangle (4, 2);
\node at (2.5, 2.5) {$R_D$};
\end{tikzpicture}}\quad
\subfloat[]{
\begin{tikzpicture}[scale = 1]
\draw (1, 2) rectangle (2, 3);
\draw (2, 2) rectangle (4, 4);
\draw (1, 3) rectangle (2, 4);
\draw (1, 1) rectangle (4, 2);
\node at (1.5, 2.5) {$R_D$};
\end{tikzpicture}}\quad
\subfloat[]{
\begin{tikzpicture}[scale = 1]
\draw (1, 1) rectangle (2, 2);
\draw (2, 1) rectangle (4, 4);
\draw (1, 2) rectangle (2, 4);
\node at (1.5, 1.5) {$R_D$};
\end{tikzpicture}}
\caption{Possible topological configurations for support subdivision in weakly singular integrals. $R_D$ is the box formed by the knots of the B-spline $\hat{B}$ enclosing the polar point $\hat{\pmb{y}}$. Depending on where such polar point is located, $R_D$ could be internal (a), on an edge (b) or in a corner (c) of $\supp(\hat{B})$. The remaining part of the support is split into rectangular regions running anti-clockwise around $R_D$ as reported in the figures.}\label{fig:supportsubdivision}
\end{figure}
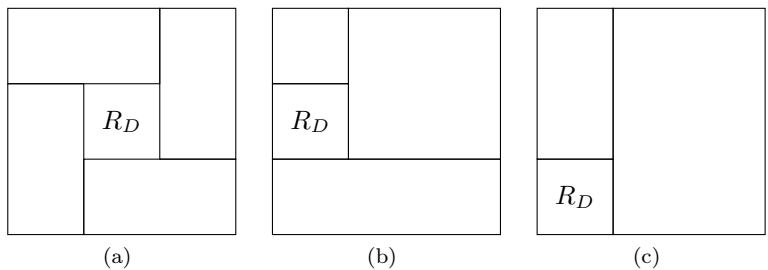

For singular integrals  $\pmb{y}$ belongs to the same $m$--th patch where $B$ is defined, as $\pmb{y} \in \supp(B)$, i.e., $\pmb{y} = \pmb{F}^{(m)}(\hat{\pmb{y}})$ and $\hat{\pmb{y}} \in \supp(\hat{B})$. Let $R_D \defeq [a_1, b_1] \times [a_2, b_2] \subseteq \supp(\hat{B})$ be the smallest rectangular region enclosing $\hat{\pmb{y}}$, 
with $\{a_1, a_2\}$  and $\{b_1, b_2\}$ being distinct knots in the local knot vectors of $\hat{B}$ in the first and second direction, respectively.  $\hat{\pmb{y}}$ may be on an edge or a corner of $R_D$, only if it lies on the boundary of $\supp(\hat{B})$. If $(d_1, d_2) \neq (0, 0)$, the remaining part of the domain $\supp(\hat{B}) \setminus R_D$ could be non-empty. This is surely true if at least one $d_k \geq 2$, for a $k \in \{1, 2\}$. We then subdivide such remaining part in a collection $\cR$ of up to four rectangular regions, running around $R_D$ in an anti-clockwise manner as in Figure \ref{fig:supportsubdivision}, depending on the position of $R_D$ in the support.

Obviously, the rectangles in $\cR$ could run clockwise around $R_D$ as well without changing essentially anything in the procedure. Moreover, such rectangles are at most $4$ independently of the the B-spline degree, leading to the computation of only up to $5$ integrals. 

Since in any rectangle $R \in \cR$ the integrand is non-singular, we can apply the procedure for non-singular integrals of Section \ref{sec:nonsing}. The only modification is in the affine transformation $\pmb{A}$, mapping $[-1, 1]^2$ into $R$ instead of $\supp(\hat{B})$.

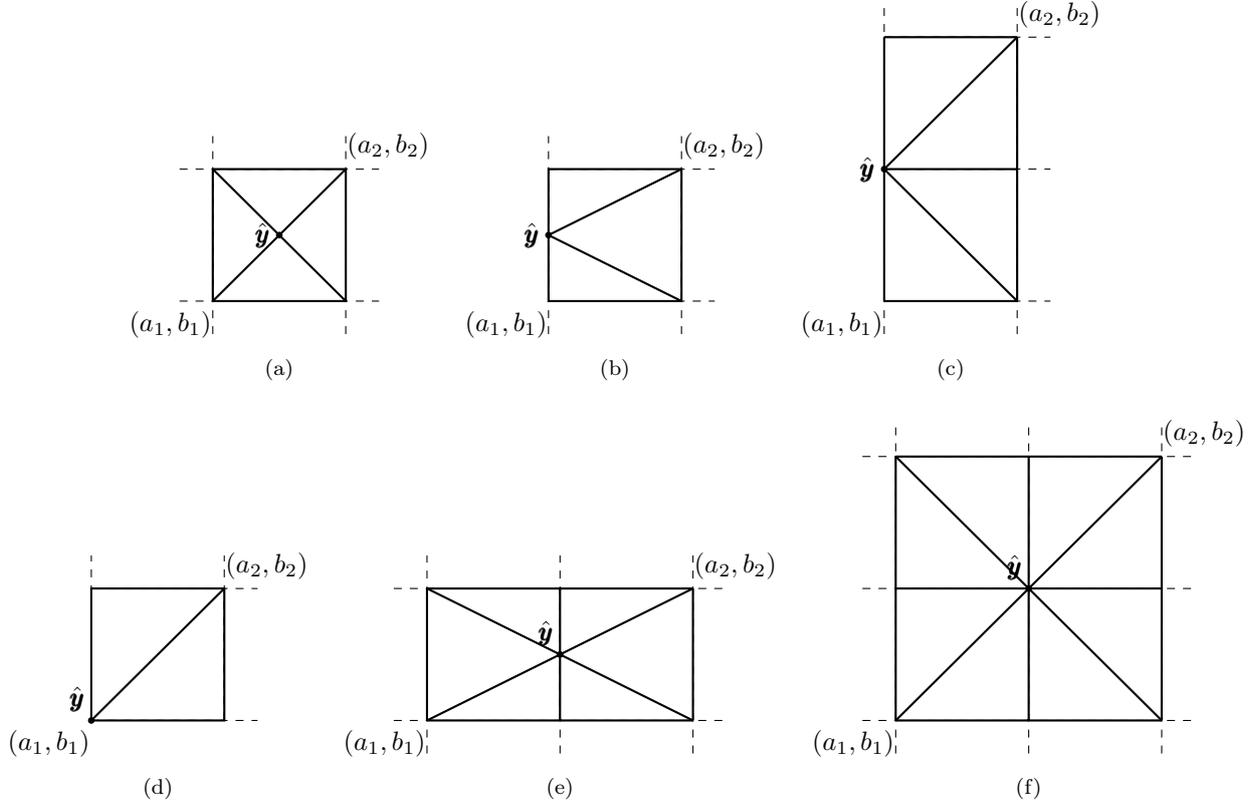
\begin{figure}
\centering
\subfloat[]{
\begin{tikzpicture}[scale = 1.75]
\draw[xstep = 1, dashed] (0, -0.25) grid (1, 1.25);
\draw[ystep = 1, dashed] (-0.25, 0) grid (1.25, 1);
\draw[thick] (0, 0) rectangle (1, 1);
\draw[thick] (0, 0) -- (1, 1);
\draw[thick] (0, 1) -- (1, 0);
\fill (0.5, 0.5) node[left]{$\hat{\pmb{y}}$} circle (0.025);
\draw[thick] (0, 0) node[below]{$(a_1, b_1)\hspace{1cm}~$};
\draw[thick] (1, 1) node[above]{$~\hspace{1cm}(a_2, b_2)$};
\end{tikzpicture}
}
\subfloat[]{
\begin{tikzpicture}[scale = 1.75]
\draw[dashed] (0, -0.25) -- (0, 1.25);
\draw[dashed] (1, -0.25) -- (1, 1.25);
\draw[ystep = 1, dashed] (0, 0) grid (1.25, 1);
\draw[thick] (0, 0) rectangle (1, 1);
\draw[thick] (0, 0.5) -- (1, 1);
\draw[thick] (0, 0.5) -- (1, 0);
\fill (0, 0.5) node[left]{$\hat{\pmb{y}}$} circle (0.025);
\draw[thick] (0, 0) node[below]{$(a_1, b_1)\hspace{1cm}~$};
\draw[thick] (1, 1) node[above]{$~\hspace{1cm}(a_2, b_2)$};
\end{tikzpicture}
}
\subfloat[]{
\begin{tikzpicture}[scale = 1.75]
\draw[step = 1, dashed] (0, -1.25) grid (1.25, 1.25);
\draw[step = 1, thick] (0, -1) grid (1, 1);
\draw[thick] (0, 0) -- (1, 1);
\draw[thick] (0, 0) -- (1, -1);
\fill (0, 0) node[left]{$\hat{\pmb{y}}$} circle (0.025);
\draw[thick] (0, -1) node[below]{$(a_1, b_1)\hspace{1cm}~$};
\draw[thick] (1, 1) node[above]{$~\hspace{1cm}(a_2, b_2)$};
\end{tikzpicture}
}
\\
\subfloat[]{
\begin{tikzpicture}[scale = 1.75]
\draw[xstep = 1, dashed] (0, 0) grid (1, 1.25);
\draw[ystep = 1, dashed] (0, 0) grid (1.25, 1);
\draw[white] (0, -0.25) -- (1, -0.25);
\draw[thick] (0, 0) rectangle (1, 1);
\draw[thick] (0, 0) -- (1, 1);
\fill (0, 0) node[above]{$\hat{\pmb{y}}\hspace{0.25cm}~$} circle (0.025);
\draw[thick] (0, 0) node[below]{$(a_1, b_1)\hspace{1cm}~$};
\draw[thick] (1, 1) node[above]{$~\hspace{1cm}(a_2, b_2)$};
\end{tikzpicture}
}
\subfloat[]{
\begin{tikzpicture}[scale = 1.75]
\draw[step = 1, dashed] (-1.25, -0.25) grid (1.25, 1.25);
\draw[step = 1, thick] (-1, 0) grid (1, 1);
\draw[thick] (-1, 0) -- (1, 1);
\draw[thick] (-1, 1) -- (1, 0);
\fill (0, 0.5) node[above]{$\hat{\pmb{y}}\hspace{0.25cm}~$} circle (0.025);
\draw[thick] (-1, 0) node[below]{$(a_1, b_1)\hspace{1cm}~$};
\draw[thick] (1, 1) node[above]{$~\hspace{1cm}(a_2, b_2)$};
\end{tikzpicture}
}
\subfloat[]{
\begin{tikzpicture}[scale = 1.75]
\draw[step = 1, dashed] (-1.25, -1.25) grid (1.25, 1.25);
\draw[step = 1, thick] (-1, -1) grid (1, 1);
\draw[thick] (-1, -1) -- (1, 1);
\draw[thick] (-1, 1) -- (1, -1);
\fill (0, 0) node[above]{$\hat{\pmb{y}}\hspace{0.25cm}~$} circle (0.025);
\draw[thick] (-1, -1) node[below]{$(a_1, b_1)\hspace{1cm}~$};
\draw[thick] (1, 1) node[above]{$~\hspace{1cm}(a_2, b_2)$};
\end{tikzpicture}
}
\caption{Possible topological configurations for subdivision in triangles of region $R_D$, according to the position of the point $\hat{\pmb{y}} \in \supp(\hat{B})$. In all the figures $R_D$ is the largest rectangular box drawn with solid lines and with $(a_1, b_1)$, $(a_2, b_2)$ as lower left and upper right corners, respectively. In (a) $\hat{\pmb{y}}$ is inside an element of the mesh in $\supp(\hat{B})$ and $R_D$ is split into 4 triangles. In (b) $\hat{\pmb{y}}$ is on the left edge of $\supp(\hat{B})$ but it is not a vertex of the mesh. $R_D$ is then split into 3 triangles. In (c) $\hat{\pmb{y}}$ is again on the left edge but this time it is also a vertex of the mesh inside  $\supp(\hat{B})$. Then $R_D$ is partitioned in 4 triangles. In (d) $\hat{\pmb{y}}$ is on a corner of $\supp(\hat{B})$ and $R_D$ can only be halved in 2 triangles. In (e) $\hat{\pmb{y}}$ is on a meshline inside $\supp(\hat{B})$ and $R_D$ is decomposed 6 triangles. Finally in (f) $\hat{\pmb{y}}$ is a mesh vertex inside $\supp(\hat{B})$ and $R_D$ is split in 8 triangles. Rotations of the splittings illustrated in figures (b)--(e) cover the remaining cases.}\label{fig:RDpartitions}   
\end{figure}

\begin{figure}
\centering
\begin{tikzpicture}
\fill (0, 0) circle (0.1);
\fill (2, -1) circle (0.1);
\fill (2, 2) circle (0.1);
\draw (0, 0)node[left]{$\hat{\pmb{y}} = (\hat{y}_1, \hat{y}_2)$} -- (2, -1)node[right]{$(a_2, b_1)$} -- (2, 2)node[right]{$(a_2, b_2)$} -- cycle;
\draw[dashed] (0, 0) --node[above, xshift = 0.35cm]{$H_T$} (2, 0);
\node at (0.5, 1) {$T$};
\end{tikzpicture}
\caption{A general representation of the triangle on the right of $\hat{\pmb{y}}$ belonging to the collection $\cT$ obtained splitting $R_D$.}
\label{fig:triangle_hat_y}
\end{figure}
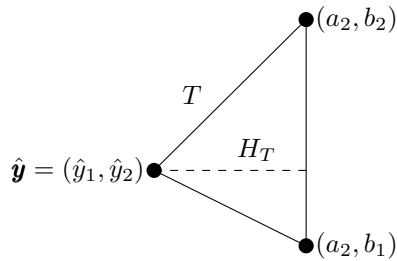

It remains to compute the integral over $R_D$. Namely, if $\hat{\pmb{x}} = (\hat{x}_1, \hat{x}_2)$ and $\hat{\pmb{y}} = (\hat{y}_1, \hat{y}_2)$, we want to compute
$$
\int_{R_D}\hat{K}_{\pmb{y}}(\hat{\pmb{x}})\hat{B}(\hat{\pmb{x}})\text{d}\hat{\pmb{x}} = \int_{a_1}^{b_1}\int_{a_2}^{b_2} \hat{K}_{\pmb{y}}(\hat{x}_1, \hat{x}_2)\hat{B}(\hat{x}_1, \hat{x}_2)\text{d}\hat{x}_1\text{d}\hat{x}_2
$$
where we have denoted by
$$
\hat{K}_{\pmb{y}}(\hat{\pmb{x}}) \defeq K(\pmb{F}^{(m)}(\hat{\pmb{x}}),  \pmb{y})J^{(m)}(\hat{\pmb{x}})
$$
for all $\hat{\pmb{x}} \in R_D$.   
In $R_D$ we shall use the \emph{Duffy desingularization} \cite{duffy}. This technique has been already considered for BEM integrals \cite{sauer, sauertesi} and it has proved to facilitate the construction of quadrature rules for singular integrals with good convergence properties. In the remaining of this subsection we summarize such procedure. Split $R_D$ into a collection $\cT$ of triangles by connecting $\hat{\pmb{y}}$ with the vertices of $R_D$, as shown in Figure \ref{fig:RDpartitions}.
Let us assume, without loss of generality, that $\hat{\pmb{y}}$ is in the interior of an element and so that the triangle $T$ on the right of $\hat{\pmb{y}}$ is in $\cT$, as depicted in Figure \ref{fig:triangle_hat_y}.
 The integration over any of the other possible triangles can be treated in a similar manner. Let $H_T$ be the height of $T$ with respect to vertex $\hat{\pmb{y}}$, as shown in Figure \ref{fig:triangle_hat_y}. The integral over $T$ is
$$
\scalebox{0.85}{$\displaystyle\iint_T \hat{K}_{\pmb{y}}(\hat{x}_1, \hat{x}_2)\hat{B}(\hat{x}_1, \hat{x}_2)\text{d}\hat{x}_1\text{d}\hat{x}_2 = \int_{\hat{y}_1}^{b_1} \int_{r_1(\hat{x}_1)}^{r_2(\hat{x}_1)} \hat{K}_{\pmb{y}}(\hat{x}_1, \hat{x}_2)\hat{B}(\hat{x}_1, \hat{x}_2)\text{d}\hat{x}_1\text{d}\hat{x}_2 \qquad \text{with }\begin{cases}
r_1(\hat{x}_1) \defeq \frac{b_1 - \hat{y}_2}{a_2 - \hat{y}_1}(\hat{x}_1 - \hat{y}_1) + \hat{y}_2\\\\
r_2(\hat{x}_1) \defeq \frac{b_2 - \hat{y}_2}{a_2 - \hat{y}_1}(\hat{x}_1 - \hat{y}_1) + \hat{y}_2.
\end{cases}$}
$$
The Duffy desingularization, in this settings, is based on the change of coordinates in the second direction $\hat{x}_2 = \hat{y}_2 + (\hat{x}_1 - \hat{y}_1)t$ so that $\text{d}\hat{x}_2 = (\hat{x}_1 - \hat{y}_1)\text{d}t$. As $\hat{x}_2$ belongs to $[r_1(\hat{x}_1), r_2(\hat{x}_1)]$ we must have $t \in [\nicefrac{(r_1(\hat{x}_1) - \hat{y}_2) }{(\hat{x}_1 - \hat{y}_1)}, \nicefrac{(r_2(\hat{x}_1) - \hat{y}_2)}{(\hat{x}_1 - \hat{y}_1)}]$. 
Such a change of variables effectively desingularizes the weakly singular kernel. Indeed, as $\hat{\pmb{x}} \to \hat{\pmb{y}}$, it holds 

$$
\hat{K}_{\pmb{y}}(\hat{x}_1,\hat{x}_2)\,\mathrm{d}\hat{x}_2 
\sim 
\frac{\mathrm{d}\hat{x}_2}{
\norm{\pmb{F}^{(m)}(\hat{\pmb{x}}) - \pmb{F}^{(m)}(\hat{\pmb{y}})}
}.
$$
Using the first-order Taylor expansion of $\pmb{F}^{(m)}$ around $\hat{\pmb{y}}$, we obtain the approximation
$$
\norm{\pmb{F}^{(m)}(\hat{\pmb{x}}) - \pmb{F}^{(m)}(\hat{\pmb{y}})}
\approx
\norm{\bbol{J}_{\pmb{F}^{(m)}}(\hat{\pmb{y}})\,(\hat{\pmb{x}} - \hat{\pmb{y}})}.
$$
Since the norm $
\norm{\hat{\pmb{x}} - \hat{\pmb{y}}}_{\bbol{J}_{\pmb{F}^{(m)}}}
\defeq
\norm{\bbol{J}_{\pmb{F}^{(m)}}(\hat{\pmb{y}})\,(\hat{\pmb{x}} - \hat{\pmb{y}})}
$ 
is equivalent to the Euclidean norm in $\mathbb{R}^3$, there exists a constant
$C_{\pmb{F}^{(m)}} > 0$ such that
$$
\norm{\bbol{J}_{\pmb{F}^{(m)}}(\hat{\pmb{y}})\,(\hat{\pmb{x}} - \hat{\pmb{y}})}
\sim
C_{\pmb{F}^{(m)}}^{-1}\, 
\norm{\hat{\pmb{x}} - \hat{\pmb{y}}}.
$$
Therefore,
\begin{equation}\label{eq:asymptotickernel}
\hat{K}_{\pmb{y}}(\hat{x}_1,\hat{x}_2)\,\mathrm{d}\hat{x}_2 
\sim
\frac{C_{\pmb{F}^{(m)}}\,\mathrm{d}\hat{x}_2}{
\sqrt{(\hat{x}_1 - \hat{y}_1)^2 + (\hat{x}_2 - \hat{y}_2)^2}
} = \frac{
C_{\pmb{F}^{(m)}}\,(\hat{x}_1 - \hat{y}_1)\,\mathrm{d}t
}{
(\hat{x}_1 - \hat{y}_1)\sqrt{1+t^2}
}
=
\frac{C_{\pmb{F}^{(m)}}\,\mathrm{d}t}{\sqrt{1+t^2}}
\end{equation}
which is a smooth function of $t$, and therefore the kernel is desingularized.
The B-spline factor, on the other hand, is polynomial when restricted to $T$, and its composition with the Duffy transformation preserves smoothness. Once the integral has been regularized, we can use a standard Gauss-Legendre quadrature rule. Hence, we consider $N_D^1$ Gauss-Legendre nodes $\{(\hat{x}_m, \hat{y}_2)\}_{m=1}^{N_D^1}$ over $H_T$. For each of such nodes $(\hat{x}_m, \hat{y}_2)$ there shall be $N_D^2$ nodes in the orthogonal direction, $\{(\hat{x}_m, t_n)\}_{m, n = 1}^{N_D^1, N_D^2}$ with $t_n$ in $[\nicefrac{r_1(\hat{x}_m)}{(\hat{x}_m - \hat{y}_1)}, \nicefrac{r_2(\hat{x}_m)}{(\hat{x}_m - \hat{y}_1)}]$. Note that, since the Gauss-Legendre nodes are all in the interior of the integration intervals, we avoid the issue of having $\hat{x}_1 = \hat{y}_1$ as quadrature node, which would have degenerated the interval where the nodes $t_n$ in the orthogonal direction are defined.

\begin{oss}
As mentioned at the beginning of this section, in order to employ the same type of quadrature formula on each B-spline support, one might consider applying the procedure introduced in Section \ref{sec:nonsing} even when the polar point lies within $\operatorname{supp}(B)$, where the integral becomes weakly singular. This could be achieved by setting $\alpha = 0$ in both coordinate directions, which corresponds to using the standard cubic Telles transformation to desingularize the integral.
However, although this approach is feasible in principle, the resulting desingularized kernel would only be continuous in the parametric domain. Consequently, the performance of our weighted interpolatory quadrature rule would be rather limited. Conversely, the Duffy transformation applied to singular integrals yields a higher degree of smoothness, as also shown in Equation \eqref{eq:asymptotickernel}. Performed numerical tests indeed confirm the superiority of this desingularization over Telles, leading to a more efficient quadrature rule.

Once the Duffy transformation is adopted for desingularization, one might wish to simplify the procedure by avoiding the aforementioned decomposition of the support into up to five rectangular regions, and instead apply Duffy directly over the entire domain. Nevertheless, the high performance of the Duffy transformation is maintained only when the additional term in the integrand, that is, the B-spline factor, is smooth over the integration domain. This condition is not satisfied, since the B-spline factor exhibits a reduced regularity along the horizontal and vertical meshlines of the parametric domain.
It should also be noted that, when the Duffy transformation is applied, the transformed B-spline is no longer a tensor product spline, and therefore cannot be treated as a weight function within our quadrature formulation.
\end{oss}
\subsubsection{Quadrature Nodes Number and Distribution}\label{sec:NS}
The approximation accuracy of the integral values should be closely tied to the fineness of the BEM discretization. To reflect this, we propose an automatic adjustment of the number of quadrature nodes used for non-singular integrals. Let $\m{N}_c = (c_1, c_2)$ denote the baseline number of quadrature nodes, given as input. We suggest increasing this number in response to the knot spacing $\m{h} = (h_1, h_2)$ in the two directions, using a double-logarithmic scaling, that is, if $h_k = 2^{-\ell_k}$ for some $\ell_k \in \NN$ and $k \in \{1, 2\}$, we define $\m{N}_{\m{h}} = (N_{h_1}, N_{h_2})$ with
$$N_{h_k} \defeq \ceil{\log_2(\ell_k)}.$$
This term $\m{N}_{\m{h}}$ will then added to $\m{N}_c$, effectively linking the quadrature density with the discretization fineness.

Additionally, since higher-degree discretizations lead to larger support regions, and therefore involve integration over larger portions of the surface $\Gamma$, we introduce a degree-scaling factor to further calibrate the quadrature nodes. Specifically, for $k \in \{1, 2\}$ we define $\m{N}_{\m{d}} = (N_{d_1}, N_{d_2})$ with
$$N_{d_k} \defeq \ceil{\log_2(d_k + 2)}.$$
The total correction becomes $\m{N}_{\m{h}} \odot \m{N}_{\m{d}}$, with $\odot$ the Hadamard (element-wise) product, which is added to $\m{N}_c$ to account for both the mesh fineness and the polynomial degree, i.e., we update
$$
\m{N}_c \leftarrow \m{N}_c + \m{N}_{\m{h}}\odot \m{N}_{\m{d}}.
$$
It is also important to note that a uniform number of quadrature nodes is not required across all integrals. While the DCT change of variables can regularize integrals enough to permit a common quadrature rule, the resulting integrand may still not be regular enough to sustain exactly the same accuracy with a uniform node distribution. To handle this, we propose adapting the number of quadrature nodes based on the geometric parameter $\alpha$, already used in the DCT transformation to measure the proximity to the singularity. The updated baseline node number $\m{c}$ is then refined into a support-specific number of quadrature points, $\m{N}_c \to \m{N}_{\supp}$, whose size depends on $\alpha$. This results in a global distribution of quadrature nodes on $\Gamma$ that is clustered near the singularity, rather than uniform.

Nevertheless, the average number of quadrature nodes per support remains low, so that the overall computational cost is not significantly increased. This is rather crucial because BEM discretizations typically require the evaluation of many integrals over all the supports partitioning $\Gamma$ and so each individual integration step must remain computationally efficient.

\section{The Stokes flow}\label{sec:stokes}
In this section, we recall the boundary integral formulation for the Stokes problem and introduce the notation adopted throughout this paper. The derivation of a Boundary Integral Equation (BIE) from the differential model is briefly described in Section \ref{sec:stokesbie}, following the approach in \cite{bie}. We also provide some integral identities that will be used in the numerical tests to evaluate the performance of our quadrature rule. Then, in Section \ref{sec:igabem}, we discuss the IgA-BEM discretization of the BIE in a multipatch surface domain setting, in a manner similar to \cite{sauer}.

The benefits of the overall methodology, which consists of reformulating the differential problem as a BIE and then discretizing it using the IgA-BEM approach, are significant. First, the use of BIE reduces the dimensionality of the problem, transforming it from a volumetric to a surface formulation. This is particularly advantageous for exterior problems with unbounded domains, as no artificial boundaries are needed for discretization. Second, volumetric meshing is often challenging for complex geometries, whereas the IgA-BEM approach leverages the fact that CAD models typically provide only the B-rep (boundary representation) of the geometry. This allows for a direct and accurate construction of the boundary mesh, aligning naturally with the integral formulation.

\subsection{Stokes differential and boundary integral equations}\label{sec:stokesbie}
Fluid flows can be classified according to their \emph{Reynolds number}, defined as the dimensionless ratio
$$
\text{Reynolds number} = \frac{\text{fluid density} \times \text{speed} \times \text{characteristic length}}{\text{viscosity}}.
$$
When the Reynolds number is small, the flow regime is referred to as \emph{slow viscous flow} or \emph{creeping flow}. However, this terminology can be misleading, as a small Reynolds number may arise not only from low flow velocity or high viscosity, but also from small characteristic length or low fluid density.

When the creeping flow of an \emph{incompressible Newtonian fluid} is however indeed sufficiently slow, the inertial forces become negligible compared to viscous forces. In these cases, the flow is then governed by the \emph{continuity equation}, which expresses the conservation of mass,
$$
\nabla \cdot \pmb{u} = 0,
$$
and the \emph{Stokes equation} (or \emph{creeping flow equation}),
$$
-\eta \Delta \pmb{u} + \nabla p = \pmb{0},
$$
where $\pmb{u}$ is the \textbf{velocity field}, $p$ is the \textbf{pressure}, and $\eta$ is the \textbf{dynamic viscosity} of the fluid. Since only the pressure gradient $\nabla p$ appears in the equations, the pressure $p$ is determined only up to an additive constant. To obtain a unique pressure field, usually a normalization condition is imposed, such as setting the mean pressure to zero
\begin{equation} \label{eq:scalp}
\int_{\Omega} p(\pmb{x}) \dx = 0\,.
\end{equation}
The Stokes equation describes a balance between viscous and pressure forces at each instant, even in the case of unsteady (time-dependent) flow. Namely, the present flow structure is independent of the flow history as it is driven only by the current configuration and by the boundary conditions on $\Gamma = \partial \Omega$. 
The most common boundary condition types are:

\begin{itemize}
    \item Dirichlet boundary condition: $\pmb{u} = \pmb{u}_D$ on $\Gamma_D \subseteq \Gamma$, where $\pmb{u}_D$ is a prescribed velocity.
    \item Neumann (traction) boundary condition: $\pmb{t} \defeq [ -p\bbol{I} + \eta\left(\bbol{J}_{\pmb{u}} + \bbol{J}_{\pmb{u}}^T\right) ]\ \m{n} = \pmb{t}_N$ on $\Gamma_N \subseteq \Gamma$, where $\m{n}$ is the outward unit normal to the boundary $\Gamma_N$, $\bbol{J}_{\pmb{u}}$ is the Jacobian matrix of $\pmb{u}$, and $\pmb{t}_N$ is a prescribed \emph{traction}.
\end{itemize}

The Dirichlet and Neumann portions of the boundary satisfy $\Gamma_D \cup \Gamma_N = \Gamma$ and $\Gamma_D \cap \Gamma_N = \varnothing$.
In order to achieve a boundary integral formulation, the knowledge of a \emph{fundamental solution} or \emph{Green's function} of such differential system is required. This is a solution of the following free space singularly forced Stokes problem:
\begin{equation}\label{eq:deltaStokes}
\left\{
\begin{array}{ll}
-\eta\Delta\pmb{u}(\pmb{x}) + \grad p(\pmb{x}) = \m{g}\delta(\pmb{x} - \pmb{y}) & \pmb{x} \in \RR^3,\\\\
\div \pmb{u}(\pmb{x}) = 0 & \pmb{x} \in \RR^3
\end{array}
\right.
\end{equation}
with $\m{g} \in \RR^3$ constant, $\pmb{y}$ an arbitrary point of the space called \emph{source point} and $\delta$ the 3D Dirac's delta. Defined $\m{r} \defeq \pmb{x} - \pmb{y}$ and $r \defeq \norm{\m{r}}$, with some mathematical effort \cite[Section 2.1]{bie}, one finds that a fundamental solution of \eqref{eq:deltaStokes} is given by
\begin{equation}\label{eq:stokesgreen}
\begin{cases}
\pmb{u}'(\pmb{x}) = \frac{1}{\eta 8\pi} \bbol{G}(\m{r})\m{g}\\\\
p'(\pmb{x}) = \frac{1}{8\pi}\m{P}(\m{r})^T\m{g} \qquad \text{with}\qquad
\displaystyle P_i(\m{r}) = 2\frac{r_i}{r^3}
\end{cases}
\end{equation}
and $\bbol{G}(\m{r})$ the $3\times 3$ symmetric matrix whose $(i,j)$--th component is
\begin{equation*}\label{eq:Stokeslet}
\bbol{G}_{ij}(\m{r}) \defeq \frac{1}{r}\delta_{ij} + \frac{r_ir_j}{r^3}.
\end{equation*}
$\bbol{G}$ is sometimes called \emph{Stokeslet} or \emph{Oseen-Burgers tensor}. 
We define the \emph{stress tensor} $\bbsigma$ related to any solution $\pmb{u}$ and $p$ of \eqref{eq:deltaStokes} as
\begin{equation*}\label{eq:stressTensor}
\bbsigma \defeq -p\bbol{I} + \eta\left(\bbol{J}_{\pmb{u}} + \bbol{J}_{\pmb{u}}^T\right).
\end{equation*}
The $(i,j)$--th entry of $\bbsigma$ is then
$$
\bbsigma_{ij} = -p\delta_{ij} + \eta(\partial_j u_i + \partial_i u_j).
$$
We note in particular that $\bbsigma$ is symmetric. Let now $\pmb{u}'$ and $p'$ be the fundamental solution \eqref{eq:stokesgreen} we found for \eqref{eq:deltaStokes},  $\bbsigma'$ the associated stress tensor and $\pmb{u}, p$ any other solution of \eqref{eq:deltaStokes}, with $\bbsigma$ the related stress tensor. Then, the following equation called \emph{Lorentz reciprocity} holds true \cite[Section 1.4]{bie} for $\pmb{x} \neq \pmb{y}$: 
\begin{equation}\label{eq:Lorentz}
\sum_{i=1}^3\sum_{j=1}^3 \partial_j(u_i'(\pmb{x})\bbsigma_{ij}(\pmb{x}) - u_i(\pmb{x})\bbsigma_{ij}'(\pmb{x})) = 0.
\end{equation}
Form the definition, we also see that the stress tensor of the fundamental solution can be written as 
$$
\bbsigma_{ij}'(\pmb{x}) = \frac{1}{8\pi}\sum_{k=1}^3\bbol{T}_{ijk}(\m{r})g_k
$$
where $\bbol{T}$ is the 3 indices tensor whose entries are
\begin{equation*}\label{eq:tensorT}
\bbol{T}_{ijk}(\m{r}) = - \delta_{ij}p_k(\m{r}) + \partial_j \bbol{G}_{ik}(\m{r}) + \partial_i \bbol{G}_{jk}(\m{r})= - 6 \frac{r_ir_jr_k}{r^5}.
\end{equation*}
We notice that $\bbol{T}$ is invariant under permutation of the indices, that is, is symmetric in any two indices. The Lorentz reciprocity \eqref{eq:Lorentz} reads, by the arbitrary of $\m{g}$, 
\begin{equation*}\label{eq:newLorentz}
\sum_{i=1}^3\sum_{j=1}^3 \partial_j\left(\bbol{G}_{ik}(\m{r})\bbsigma_{ij}(\pmb{x}) - \eta u_i(\pmb{x}) \bbol{T}_{ijk}(\m{r})\right) = 0 \qquad \forall\,k.
\end{equation*}
Let now $\Omega\subseteq \RR^3$ be an open domain with Lyapunov boundary $\partial\Omega = \Gamma$ and let $\pmb{y} \in \Omega$. Let also $B_\eps(\pmb{y})$ be the ball centered at $\pmb{y}$ of radius $\eps>0$. When integrating over $\Omega \setminus B_\eps(\pmb{y})$, by the divergence theorem we have
{\footnotesize\begin{equation*}
\sum_{i=1}^3 \int_{\Omega} \sum_{j=1}^3 \partial_j\left(\bbol{G}_{ik}(\m{r})\bbsigma_{ij}(\pmb{x}) - \eta u_i(\pmb{x}) \bbol{T}_{ijk}(\m{r})\right)\dx = \sum_{i=1}^3 \int_{\Gamma\cup \partial B_\eps(\pmb{y})} \sum_{j=1}^3 n_j(\pmb{x})\left(\bbol{G}_{ik}(\m{r})\bbsigma_{ij}(\pmb{x}) - \eta u_i(\pmb{x}) \bbol{T}_{ijk}(\m{r})\right)\dsigma
\end{equation*}}
for all $k$, where $\m{n}(\pmb{x})$ is the normal inward $\Omega$ at $\pmb{x} \in \Gamma \cup \partial B_\eps(\pmb{y})$. As we are far from $\pmb{y}$, the Lorentz reciprocity still holds true so that the left hand-side is zero. Hence we get
\begin{equation*}
\sum_{i=1}^3\sum_{j=1}^3 \int_{\Gamma\cup \partial B_\eps(\pmb{y})}  n_j(\pmb{x})\left(\bbol{G}_{ik}(\m{r})\bbsigma_{ij}(\pmb{x}) - \eta u_i(\pmb{x}) \bbol{T}_{ijk}(\m{r})\right)\dsigma = 0 \qquad \forall\, k.
\end{equation*}
By splitting the integration in $\Gamma$ and $\partial B_\eps(\pmb{y})$ and letting $\eps \to 0$ to eliminate the artificial latter boundary around $\pmb{y}$, with a change of variables to the solid angle in the integral over such infinitesimal ball, one ultimately shows \cite[Section 2.3]{bie} that for any $k$
\begin{equation}\label{eq:birfcomponent}
u_k(\pmb{y}) = -\frac{1}{8\pi \eta} \sum_{i=1}^3\sum_{j=1}^3 \int_{\Gamma} n_j(\pmb{x})\bbol{G}_{ik}(\m{r})\bbsigma_{ij}(\pmb{x})\dsigma + \frac{1}{8\pi}\sum_{i=1}^3\sum_{j=1}^3\int_{\Gamma} n_j(\pmb{x})u_i(\pmb{x}) \bbol{T}_{ijk}(\m{r})\dsigma.
\end{equation}
If we define, for any $\pmb{x} \in \Gamma$, the scalar functions $U(r) \defeq \nicefrac{1}{r}$ and $H(\m{r}; \m{n}(\pmb{x})) \defeq \nicefrac{\langle\m{r}, \m{n}(\pmb{x})\rangle}{r^3}$, the dyadic matrix $\bbol{r} \defeq \nicefrac{\m{r}\m{r}^T}{r^2}$ and the \textbf{surface force}, also called \textbf{traction}, vector $\pmb{t}(\pmb{x}) \defeq \bbsigma(\pmb{x})\m{n}(\pmb{x})$, we notice that
$$
\bbol{G}(\m{r}) = \left(\bbol{I} + \bbol{r}\right) \qquad\text{and}\qquad \sum_{j=1}^3 n_j(\pmb{x}) \bbol{T}_{ijk}(\m{r}) = - 6H(\m{r}; \m{n}(\pmb{x}))\bbol{r}
$$
and by using the symmetric properties of $\bbol{G}$ and $\bbol{T}$ to commute the matrix-vector products, we can rewrite Equation \eqref{eq:birfcomponent} in a vector form as
\begin{equation}\label{eq:birf}
\pmb{u}(\pmb{y}) = - \frac{1}{8\pi\eta}\int_\Gamma U(r)\left(\bbol{I} + \bbol{r}\right)\pmb{t}(\pmb{x})\dsigma - \frac{3}{4\pi} \int_\Gamma H(\m{r}, \m{n}(\pmb{x})) \bbol{r}\pmb{u}(\pmb{x})\dsigma.
\end{equation}
This equation provides the value of the velocity field $\pmb{u}$ of the flow for all $\pmb{y}$ in the open domain $\Omega$ in terms of two boundary integrals, called \textbf{single- and double-layer potentials}, respectively. Whereas, the matrix functions $U(r)\left(\bbol{I} + \bbol{r}\right)$ and $H(\m{r}, \m{n}(\pmb{x})) \bbol{r}$ are referred to as \textbf{single- and double-layer kernels}, respectively. Therefore, Equation \eqref{eq:birf} is called \textbf{boundary integral representation formula}. Then, when we let $\pmb{y}$ approach $\Gamma$ the double-layer potential verifies \cite[Equation 2.3.12]{bie}
$$
\lim_{\pmb{y} \to \Gamma}  \int_\Gamma H(\m{r}, \m{n}(\pmb{x})) \bbol{r}\pmb{u}(\pmb{x})\dsigma = -\frac{2\pi}{3} \pmb{u}(\pmb{y}) + \int_\Gamma H(\m{r}, \m{n}(\pmb{x})) \bbol{r}\pmb{u}(\pmb{x})\dsigma.
$$ 
Replacing the double-layer with such limit in Equation \eqref{eq:birf} we get
\begin{equation}\label{eq:BIEmatrixform}
\pmb{u}(\pmb{y}) = - \frac{1}{4\pi\eta}\int_\Gamma U(r)\left(\bbol{I} + \bbol{r}\right)\pmb{t}(\pmb{x})\dsigma - \frac{3}{2\pi} \int_\Gamma H(\m{r}, \m{n}(\pmb{x})) \bbol{r}\pmb{u}(\pmb{x})\dsigma
\end{equation}
which is the implicit representation formula for $\pmb{y} \in \Gamma$, called \textbf{Boundary Integral Equation} (BIE). We observe that, under our assumption that the boundary is a Lyapunov surface, the entries of the matrix integrands of the BIE are weakly singular functions, according to the definitions of $U, H$ and $\bbol{r}$. If we solve for $\pmb{u}$ such implicit formula, then we know the value of the velocity field for any point $\pmb{x} \in \Gamma$. With this information, we have all the data to apply the boundary representation formula \eqref{eq:birf} to determine the velocity field $\pmb{u}$ at any point $\pmb{y}$ in the domain $\Omega$, just by integrating over $\Gamma$ the layers. The pressure at any point $\pmb{y} \in \Omega$ is also readily determined from the knowledge of $\pmb{u}$ and $\pmb{t}$ on $\Gamma$, by differentiating the boundary integral representation formula for $\pmb{u}$ and using the Stokes equation on the fundamental solution, for $\pmb{x} \neq \pmb{y}$. This leads to the \textbf{boundary representation formula for the pressure}
\begin{equation*}
\begin{split}
p(\pmb{x}) &= -\frac{1}{8 \pi}\ \int_{\Gamma}
{\bf P}(\m{r})^T \pmb{t}(\pmb{x}) \dsigma +\frac{\eta}{8 \pi} \ \int_{\Gamma} \pmb{u}^T(\pmb{x}) \bbol{P}\m{n}(\pmb{x})\dsigma + C\\
&= -\frac{1}{8 \pi}\ \int_{\Gamma}
{\bf P}(\m{r})^T \pmb{t}(\pmb{x}) \dsigma +\frac{\eta}{8 \pi} \ \int_{\Gamma} \pmb{u}^T(\pmb{x}) 4(U(r))^3(3\bbol{r} - \bbol{I})\m{n}(\pmb{x})\dsigma + C
\end{split}
\end{equation*}
with $C \in \RR$ a constant, which will be set to zero when imposing the normalization requirement in \eqref{eq:scalp}, and $\bbol{P}$ the following $2$ indices tensor,
$$\bbol{P}_{ij} = \frac{4}{r^3} \left( \frac{3r_ir_j}{r^2} - \delta_{ij} \right).$$
We finally underline that the Neumann and Dirichlet boundary conditions, provided in the formulation of the Stokes problem, enter the BIE \eqref{eq:BIEmatrixform} in the single- and double-layer integrals yielding $\pmb{t}$ or $\pmb{u}$ over $\Gamma$, respectively, and hence the value of one of the two integrals (and the left hand-side in case of Dirichlet boundary conditions). For mixed boundary conditions, one can  split $\Gamma$ into pieces where Neumann and Dirichlet boundary conditions are imposed, $\Gamma = \Gamma_N \cup \Gamma_D$, and solve the BIE for $\pmb{y} \in \Gamma_N$ and $\pmb{y} \in \Gamma_D$ simultaneously. 

Besides the BIE of the Stokes problem we also introduce the following two identities \cite[Equations (2.1.4) and (2.1.12)]{bie}, which will be later used to numerically test the approximation power and robustness of the quadrature rule proposed in Section \ref{sec:quadrule}.  The first identity can be easily obtained from the continuity equation and the divergence theorem, applied to the velocity field of the fundamental solution \eqref{eq:stokesgreen}:
\begin{equation*}\label{eq:identity1}
\int_\Gamma \bbol{G}(\m{r})\m{n}(\pmb{x})\dsigma = \pmb{0}\qquad\forall\, \pmb{y} \in \Gamma.
\end{equation*}
The second identity is also a consequence of the divergence theorem. By plugging the fundamental solution into the Stokes equation and recalling the definition of $\bbol{T}$, we get
\begin{equation*}\label{eq:identity2}
-\frac{1}{4\pi} \int_\Gamma \sum_{k=1}^3 \bbol{T}_{ijk}n_k(\pmb{x})\dsigma = \bbol{I} \qquad\forall\, \pmb{y}\in \Gamma.
\end{equation*}
We can rewrite these two identities as
\begin{equation}\label{eq:identities}
\int_\Gamma (\bbol{I} +\bbol{r})\m{n}(\pmb{x})\dsigma = \pmb{0};\qquad -\frac{1}{4\pi}\int_{\Gamma} H(\m{r}, \m{n}(\pmb{x}))\bbol{r}\dsigma=\bbol{I} \qquad\forall\, \pmb{y} \in \Gamma.
\end{equation}
We underline that the functions in these integrals are weakly singular with respect to $\pmb{y} \in \Gamma$.  

\begin{oss}
If the boundary is less regular than Lyapunov, the unit normal may be discontinuous or undefined along edges, corners, or other non-smooth portions of $\Gamma$. 
As a consequence, the double-layer potential is not well defined pointwise everywhere on the boundary. 
To handle this situation, a broader analytical framework based on Sobolev trace spaces is required, together with additional numerical care near geometric singularities.
\end{oss}

\subsection{Isogeometric Discretization via Collocation}\label{sec:igabem}
We now describe the discretization process via collocation of the Stokes BIE \eqref{eq:BIEmatrixform}, as well as of the identities \eqref{eq:identities}, in the IgA-BEM setting. Therefore, we assume hereafter that $\Gamma$ is a multipatch NURBS geometry composed of $M$ quadrilateral patches, $\Gamma = \cup_{m=1}^M \Gamma^{(m)},$ with $ \mathring{\Gamma}^{(\ell)}   \cap  \mathring{\Gamma}^{(m)} = \varnothing$ if $\ell \neq m,$ where $\partial \Gamma^{(\ell)} \cap \partial \Gamma^{(m)}$ is empty or it is a common corner or a common boundary edge.  
\\
For $m = 1, \ldots, M,$ let $\m{B}^{(m)}\defeq (B_{(ij)}^{(m)})_{(ij) = 1}^{N^{(m)}}$ be the column vector of the basis functions defined in $\Gamma^{(m)}$ generating the $m$--th pushforward space associated to a bivariate B-spline space on the reference domain $[0\,,\,1]^2$, with $(ij)$ the flattened index of the pair of indices $(i, j)$. 
Furthermore, let $\pmb{y}_{(ij)}^{(m)} \,, (ij)=1,\ldots,N^{(m)}\,,$ denote the corresponding $N^{(m)}$ points on $\Gamma^{(m)}$ obtained as pushforward of the tensor product Greville points of such $m$--th parametric bivariate B-spline space.  
Denoting as $\bbol{I}$ the $3\times 3$ identity matrix, we define the following matrix of size $3 \times 3N^{(m)}$,
\begin{equation*}
\bbol{B}^{(m)} \defeq (\m{B}^{(m)})^T \otimes \bbol{I} = \left[
\begin{array}{ccc|ccc|c|ccc}
B_1^{(m)}	&		&		& B_2^{(m)}	 &		 &		 &		 & B_{N^{(m)}}^{(m)}	 &		    & \\
	& B_1^{(m)}	&		& 		 & B_2^{(m)}	 &		 & \cdots &			 & B_{N^{(m)}}^{(m)} & \\
	&		& B_1^{(m)}	&		 &		 & B_2^{(m)}	 &		 &			 &		    & B_{N^{(m)}}^{(m)}		
\end{array}
\right].
\end{equation*}
With this notation, let us first introduce the discretization of the BIE \eqref{eq:BIEmatrixform} assuming only Dirichlet boundary conditions for the original Stokes problem. Then, the velocity field $\pmb{u}$ is known on $\Gamma$ and we only need to approximate there the traction $\pmb{t}$ using the following formula,
$$
\pmb{t}(\pmb{x}) \approx \pmb{t}^h(\pmb{x}) \defeq \bbol{B}^{(m)}(\pmb{x})\pmb{\uptau}^{(m)}, \text{ when }\pmb{x} \in \Gamma^{(m)}\,, \qquad\text{with } \pmb{\uptau}^{(m)} = \begin{bmatrix}
\vdots \\ \uptau_{\ell}^{1, (m)} \\ \uptau_{\ell}^{2, (m)} \\ \uptau_{\ell}^{3, (m)} \\ \vdots
\end{bmatrix} \text{ for }\ell=1, \ldots, N^{(m)}
$$
 Namely, the $i$--th component of the traction field $\pmb{t}$, denoted by $t^i$, is approximated on $\Gamma^{(m)}$ as the following linear combination
$$
t^i(\pmb{x}) \approx \sum_{\ell = 1}^{N^{(m)}} \uptau_\ell^{i, (m)} B_\ell^{(m)}(\pmb{x})\,, \qquad  \pmb{x} \in \Gamma^{(m)}\,, \qquad i=1,2,3.
$$
Let then $\pmb{\uptau} \defeq (\pmb{\uptau}^{(m)})_{m=1}^M$ be the vector collecting all the degrees of freedom (DOFs), obtained by stacking the column vectors $\pmb{\uptau}^{(m)}$. The collocation method consists in imposing the fulfillment of the BIE \eqref{eq:BIEmatrixform} pointwise in as many points as the dimension of the approximation space in order to obtain a squared linear system. The standard approach in the IgA collocation approach is to use the Greville physical points $\pmb{y}_{(ij)}^{(m)} \,, (ij)=1,\ldots,N^{(m)}\,, m=1,\ldots,M$, which lead to the following equation:
\begin{equation}\label{eq:dirichlet}
\m{C} = \bbol{S}\pmb{\uptau} + \m{D}.
\end{equation}

Setting $N \defeq \sum_{m = 1}^M N^{(m)}\,,$ $\m{C}$ is the vector with $3N$ components obtained by evaluating $\pmb{u}$ at all the $\pmb{y}_{(ij)}^{(m)}$ for all $(ij) = 1, \ldots, N^{(m)}\,, m = 1, \ldots, M$ and $\bbol{S}$ is a square matrix of size $3N \times 3N$, collocating the single-layer integral. 
$\bbol{S}$ can be decomposed in $M^2$ blocks, each one corresponding to a pair $(m_1, m_2), m_1,m_2 =1,\ldots,M,$ of patches. The block $\bbol{S}^{(m_1, m_2)}$ related to the pair $(m_1,m_2)$ has size $3N^{(m_1)}\times 3N^{(m_2)}\,,$ with $m_1$ denoting the patch where the involved collocation points are located and $m_2$ to the patch of the considered degrees of freedom. Each block can in turn be partitioned into $3\times 3$ blocks corresponding to the $3 \times 3$ matrix integrand for a fixed pushforward Greville point in patch $m_1$ and a fixed pushforward B-spline in the approximation space on patch $m_2$. More specifically, defined for all $m =1, \ldots, M$ the 3-index $\m{(ij)} \defeq (3(ij) - 2, 3(ij) - 1, 3(ij))^T$ for $(ij)$ in the the range $1, \ldots, N^{(m)}$ and fixed $(m_1, m_2)$, the $3\times3$ block $\bbol{S}_{\m{(i_1j_1)}, \m{(i_2j_2)}}^{(m_1, m_2)}$ corresponding to the 3-indices $\m{(i_1j_1)} $ and $\m{(i_2j_2)}$ is given by
$$
\bbol{S}_{\m{(i_1j_1)}, \m{(i_2j_2)}}^{(m_1, m_2)} = \frac{1}{4\pi\eta} \int_{\supp(B)_{(i_2j_2)}^{(m_2)}} U(r) (\bbol{I} + \bbol{r})B_{(i_2j_2)}^{(m_2)}(\pmb{x})\dsigma
$$
with $r$ and $\bbol{r}$ computed with respect to $\pmb{y}_{(i_1j_1)}^{(m_1)}$. Therefore, the collection of 3 rows corresponding to the 3-index $\m{(i_1j_1)}$ in the block $(m_1, m_2)$ shall be
$$
\bbol{S}_{\m{(i_1j_1)}, :}^{(m_1, m_2)} = \frac{1}{4\pi\eta} \int_{\supp(B)_{(i_2j_2)}^{(m_2)}} U(r) (\bbol{I} + \bbol{r})\bbol{B}^{(m_2)}(\pmb{x})\dsigma
$$
where we have used the Python/MATLAB ``colon'' notation to indicate that we are running over all the $\{\m{(i_2j_2)}\,:\, (i_2j_2) = 1, \ldots, N^{(m_2)}\}$. Its entries are computed using the quadrature rule proposed in Section \ref{sec:quadrule} after the pullback operation. 
Finally, the vector $\m{D}$ in \eqref{eq:dirichlet} is the collocation of the double-layer integral of the BIE \eqref{eq:BIEmatrixform} at all the Greville points $\pmb{y}_{(ij)}^{(m)}$ for all $(ij) = 1, \ldots, N^{(m)}$ and all $m = 1, \ldots, M$. Hence, here as well we have to approximate the value of integrals, although the function to integrate is known because the boundary velocity field $\pmb{u}|_{\Gamma}$ is given as datum. The accuracy used for the computation of the entries of $\m{D}$ should go hand in hand with the discretization fineness. Thanks to the partition of unity of B-splines, we can use also in these instances the quadrature rule of Section \ref{sec:quadrule}. In fact, we simply have
\begin{equation*}
\begin{split}
- \frac{3}{4\pi} \int_{\Gamma^{(m)}} H(\m{r}, \m{n}(\pmb{x})) \bbol{r}\pmb{u}(\pmb{x})\dsigma &= - \frac{3}{4\pi} \sum_{(ij) = 1}^{\mathcal{N}^{(m)}}\int_{\Gamma^{(m)}} H(\m{r}, \m{n}(\pmb{x})) \bbol{r}\pmb{u}(\pmb{x})\mathcal{B}_{(ij)}^{(m)}(\pmb{x})\dsigma \\&= - \frac{3}{4\pi} \sum_{(ij) = 1}^{\mathcal{N}^{(m)}}\int_{\supp(\mathcal{B}_{(ij)}^{(m)})} H(\m{r}, \m{n}(\pmb{x})) \bbol{r}\pmb{u}(\pmb{x})\mathcal{B}_{(ij)}^{(m)}(\pmb{x})\dsigma
\end{split}
\end{equation*}
with $\pmb{\mathcal{B}}^{(m)} = (\mathcal{B}_{(ij)}^{(m)})_{(ij)=1}^{\mathcal{N}^{(m)}}$ any pushforward B-spline basis on the $m$--th patch, as the partition of unity property holds regardless of the specific set chosen. In particular, $\pmb{\mathcal{B}}^{(m)}$ can be completely unrelated to $\m{B}^{(m)}$, meaning that they may have different degrees and they may be defined as pushforward of B-splines on different meshes. The partition of unity not only allows us to use the quadrature rule of Section \ref{sec:quadrule}, but it also partitions the domain of integration and thereby it facilitates the achievement of a suitable accuracy. Thus, vector $\m{D}$ of Equation \eqref{eq:dirichlet} is obtained by summing the columns of a matrix $\bar{\bbol{D}}$, i.e.,
$\m{D} = \bar{\bbol{D}}\m{1},$
with $\bar{\bbol{D}}$ that can be divided into $3\times 1$ blocks of the type
$$
\bar{\bbol{D}}_{\m{(i_1j_1)}, (i_2j_2)}^{(m_1, m_2)} = - \frac{3}{4\pi} \int_{\supp(\mathcal{B}_{(i_2j_2)}^{(m_2)})} H(\m{r}, \m{n}(\pmb{x})) \bbol{r}\pmb{u}(\pmb{x})\mathcal{B}_{(i_2j_2)}^{(m_2)}(\pmb{x})\dsigma
$$
where $r, \bbol{r}$ are computed with respect to the physical Greville points $(\pmb{y}_{(i_1j_1)}^{(m_1)})_{(i_1j_1) = 1}^{N^{(m_1)}}$ of the appoximant space $\m{B}^{(m_1)}$. We underline that summing over the columns of $\bar{\bbol{D}}$ corresponds in summing over the indices $(i_2j_2) \in \{1, \ldots, \mathcal{N}^{(m_2)}\}$ and all $m_2 \in \{1, \ldots, M\}$.
This concludes the discretization of the BIE \eqref{eq:BIEmatrixform} when Dirichlet conditions are provided. 
We remark that we chose a $C^{-1}$ discretization, that is, the obtained discrete approximation $\pmb{t}^h$ of the traction is discontinuous, with finite jumps, along the patch interfaces, since no constraints are adopted to ensure global continuity. This greatly simplifies the implementation on non-conforming meshes, since no gluing constraints are required for $C^0$ smoothness across patches with meshes of different resolutions. In the BEM setting, such discontinuities are usually harmless: the jumps occur only along patch edges, which have zero measure in the surface integrals, and therefore have a negligible impact on accuracy. The $C^{-1}$ choice thus trades a marginal loss of regularity for greater flexibility and simplicity.

The integrals in the identities \eqref{eq:identities} can be computed following the same procedure used for computing vector $\m{D}$ in \eqref{eq:dirichlet}, that is, by plugging in a B-spline factor in the integrand and splitting the integral in a sum, thanks to the partition of unity property, then applying the quadrature rule of Section \ref{sec:quadrule} to each addend and summing up the columns of the matrix obtained from the subdivision of $\Gamma$ over  the B-spline supports and a corresponding set of Greville points.

In case of Neumann boundary condition, the discretization is the mirrored version of the procedure described for the Dirichlet problem. Namely, this time the traction $\pmb{t}$ is known on $\Gamma$. The collocation method applied to the BIE \eqref{eq:BIEmatrixform} leads to
\begin{equation}\label{eq:neumann}
\bbol{C}\pmb{\upupsilon} = \m{S} + \bbol{D}\pmb{\upupsilon}.
\end{equation}
with $\bbol{C}$ and $\m{S}$ the collocation matrix and the evaluation vector of the single-layer integral of the BIE \eqref{eq:BIEmatrixform}, respectively, at all the physical Greville points over all the pushforward B-splines of all the patches. $\pmb{\upupsilon}$ is the unknown vector of B-spline coefficients necessary to approximate the velocity field $\pmb{u}$ in $\Gamma$ and $\bbol{D}$ is the matrix obtained by collocating the double-layer of \eqref{eq:BIEmatrixform}.

Finally, for mixed boundary conditions, we evaluate the Dirichlet boundary equation \eqref{eq:dirichlet} at the collocation points on $\Gamma_D$ and the Neumann boundary equation \eqref{eq:neumann} at the collocation points on $\Gamma_N$. The corresponding degrees of freedom in the resulting linear system are therefore interlaced and provide simultaneously the approximations $\pmb{t}^h$ on $\Gamma_D$ and $\pmb{u}^h$ on $\Gamma_N$, i.e., the missing boundary data.

As final remark, we emphasize that the collocation method requires the collocation points to be all distinct in order to avoid repeated rows in the discretization systems \eqref{eq:dirichlet}--\eqref{eq:neumann}. The Greville points lying on the interfaces of patches with conforming meshes violate this condition. Hence, for these instances, we shall replace the standard Greville points with what is called improved Greville points \cite{improvedgreville}, that is, a shift of the boundary Greville points in the interior of each patch, in order to have all distinct rows in the system matrices. In particular, this modification is also needed to have the block diagonal structure of matrix $\bbol{C}$. 

\section{Numerical examples}\label{sec:numericaltests}
Here we test the quadrature rule proposed in Section \ref{sec:quadrule} for the verification of the identities \eqref{eq:identities}, related to the kernels of the Stokes problem, as well as for checking their influence in the discretization of the BIE \eqref{eq:BIEmatrixform} with the IgA-BEM collocation. In all tests, the function $\alpha$ used for the dependence of the quadrature rule on the distance between the integration domain and the polar point is taken from the family \eqref{eq:alphaclass} with $\gamma = 8$. The number of Chebyshev nodes is set according to the B-spline degree $\m{d} = (d_1, d_2)$ as $\m{N}_c = 2 \m{d} + \pmb{3}$. This dependence on $\m{d}$ is done to make the quadrature setting self-regulating once the discretization is fixed and to take into account the enlargement of the supports, i.e., of the domains of integration, with the discretization degree. We underline that in particular, for $\m{d} = (0, 0)$ we get $\m{N}_c = (3, 3)$, which is therefore the minimal number considered. As described in Section \ref{sec:NS}, the algorithm then adjusts the input $\m{N}_c$ to a support-specific number of quadrature points $\m{N}_{\supp}$, for the non-singular integrals. In order to have a fully automatic procedure, for the singular integrals we follow the rule of thumb used already in \cite{sauer} and set $\m{N}_D = (N_D^1, N_D^2) = 2\m{N}_c$. Finally, for simplicity, and because anisotropic discretizations are unnecessary in the tests considered, we highlight that we choose to use the same polynomial degree $\mathbf{d} = (d, d)$ and the same discretization step $\mathbf{h} = (h, h)$ in both directions and in all patches throughout the numerical experiments.
\subsection{Tests on Stokes single- and double-layer kernel identities}
Fixed a collocation point $\pmb{y}_{(ij)}^{(m)}$, lying in patch $\Gamma^{(m)}$ for some $m \in \{1, \ldots, M\}$, let us indicate as $\m{I}_{\pmb{y}_{(ij)}^{(m)}}^{\sl}, \bbol{I}_{\pmb{y}_{(ij)}^{(m)}}^{\dl}$ the value computed with our quadrature rule, following the procedure explained in Section \ref{sec:igabem}, for the single- and double-layer identities, respectively, reported in \eqref{eq:identities} with respect to $\pmb{y}_{(ij)}^{(m)}$.  By recalling that the exact value of the two integrals are $\pmb{0}$ and $\bbol{I}$ respectively, regardless of $\pmb{y}_{(ij)}^{(m)}$, we define the quadrature error for collocation point $\pmb{y}_{(ij)}^{(m)}$ as
$$
e_{\pmb{y}_{(ij)}^{(m)}}^{\sl} \defeq \norm{\m{I}_{\pmb{y}_{(ij)}^{(m)}}^{\sl}}_{\text{F}}\qquad e_{\pmb{y}_{(ij)}^{(m)}}^{\dl} \defeq \norm{\bbol{I}_{\pmb{y}_{(ij)}^{(m)}}^{\dl} - \bbol{I}}_{\text{F}}
$$
with $\norm{\cdot}_{\text{F}}$ the Frobenius norm. Then, we define the mean quadrature errors for the the two identites as the average over the collocation points, i.e.,
$$
e^{\sl} \defeq \frac{1}{N} \sum_{m=1}^{M} \sum_{(ij)=1}^{N^{(m)}} e_{\pmb{y}_{(ij)}^{(m)}}^{\sl}, \qquad e^{\dl} \defeq \frac{1}{N} \sum_{m=1}^{M} \sum_{(ij)=1}^{N^{(m)}} e_{\pmb{y}_{(ij)}^{(m)}}^{\dl}
$$
where we recall that $N \defeq \sum_{m = 1}^M N^{(m)}$. Our quadrature rule are tested by measuring $e^{\sl}$ and $e^{\dl}$ against the mean number $\bar{n}_{qp}$ of quadrature points per support, that is, the total number of distinct quadrature points  divided by the number of computed integrals, which is the number of supports over $\Gamma$. We recall that, for our rule, the quadrature points are given by the Chebyshev nodes used for the polynomial interpolation of the kernels in the non-singular regions, not containing the collocation point, and the Gauss-Legendre nodes in the Duffy triangles for the singular region on $\Gamma$. 

\begin{figure}
\centering
\subfloat[]{\includegraphics[width = 0.3\textwidth, page = 1]{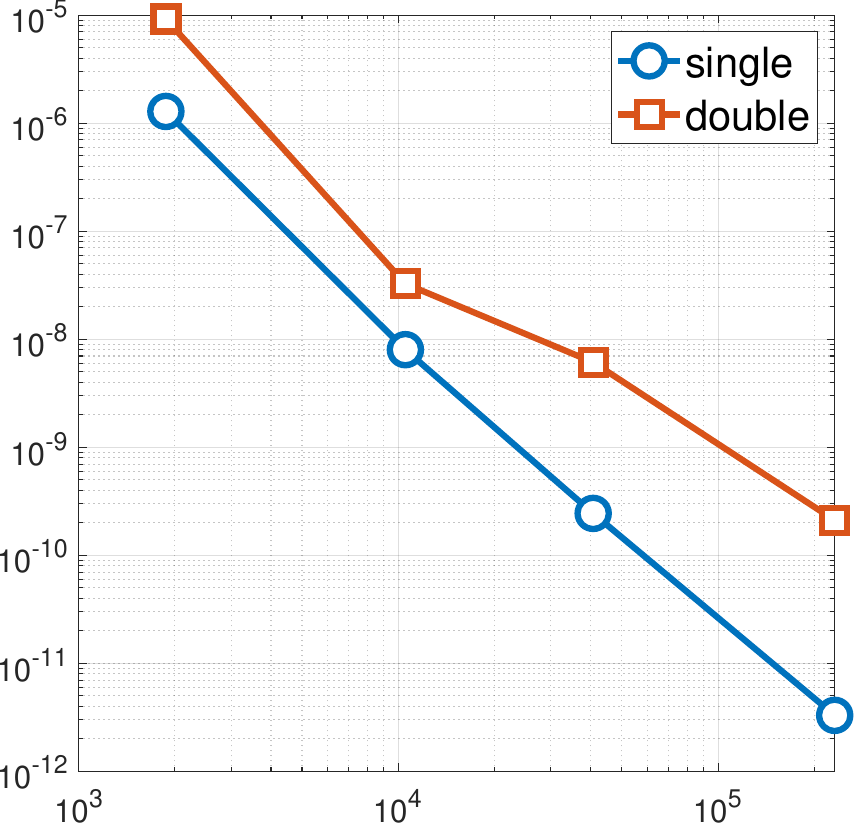}}\quad
\subfloat[]{\includegraphics[width = 0.3\textwidth, page = 2]{tests}}\quad
\subfloat[]{\includegraphics[width = 0.3\textwidth, page = 3]{tests}}
\caption{Performance of the proposed quadrature rule for the Stokes kenerl identities \eqref{eq:identities} on a sphere. The sphere is described as a 6-patches quartic NURBS. In (a) we consider B-spline sets of degree $\m{d} = (0, 0)$ on each patch and the mesh resolution of $2^\ell\times 2^\ell$ elements is taken for $\ell = 2, \ldots, 5$. On the $x$-axis the mean number of quadrature points per collocation point $\bar{n}_{\text{qp}}$ and on the $y$-axis the values of $e^{\sl}$ (round markers) and $e^{\dl}$ (square markers). In (b) we compare the values of $e^{\sl}$ for different B-spline degrees, namely $\m{d} = \pmb{0}, \ldots, \pmb{3}$, to check the robustness of the quadrature in $\m{d}$. Here on the $x$-axis we have the mesh fineness, i.e., $\ell =2, \ldots, 5$. In (c) we check the robustness in $\m{d}$ for the double-layer kernel as well, i.e., the values of $e^{\dl}$. In (b)--(c) we can notice a similar decay regardless of $\m{d}$.}\label{fig:id}
\end{figure}

In the experiments reported in this subsection $\Gamma$ is always a sphere written as a 6-patches NURBS of rational degree 4, since this is probably the standard exact representation of a sphere \cite[Section 6]{tilingsphere}, even if alternatives can be adopted.
In Figure \ref{fig:id} (a) we show the results for such $\Gamma$ and for B-spline bases of degree zero over the patches, defined on meshes of increasing resolutions. 
For these settings, the integration is done element-wise as the B-spline supports coincide with the mesh elements. This allows also for a comparison with \cite[Figure 20 (a)--(b)]{sauer}, which is the state-of-the-art for this problem to the best of our knowledge: it can be clearly seen a much faster decay of the errors for a comparable mean number of quadrature nodes. In Figure \ref{fig:id} (b)--(c), we instead test the robustness with respect to higher degrees of the B-spline bases. 

\subsection{Simulations of the exterior Stokes flow caused by a rotating sphere}\label{sec:rotsphere}
Assume to have a solid sphere $S$ of radius $R$, with surface $\Gamma = \partial S$, immersed in a Newtonian fluid with dynamic viscosity $\eta$. The sphere is rotating around the $z$-axis with an angular velocity module $\omega$, causing the motion of the surrounding fluid and hence the creation of a creeping flow, orthogonal to the $z$-axis. 

\begin{figure}
\centering
\subfloat[]{\includegraphics[width = 0.22\textwidth, page = 6]{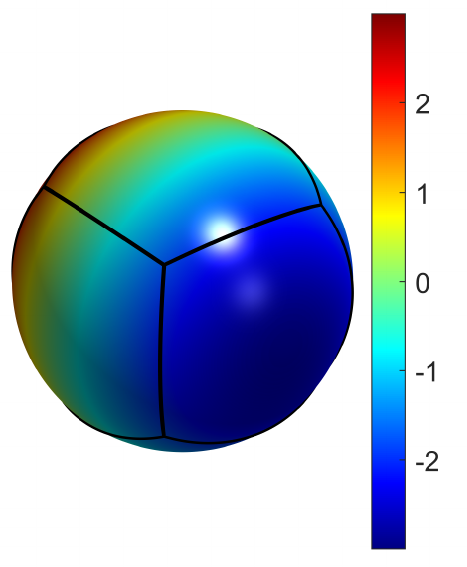}}\quad
\subfloat[]{\includegraphics[width = 0.22\textwidth, page = 7]{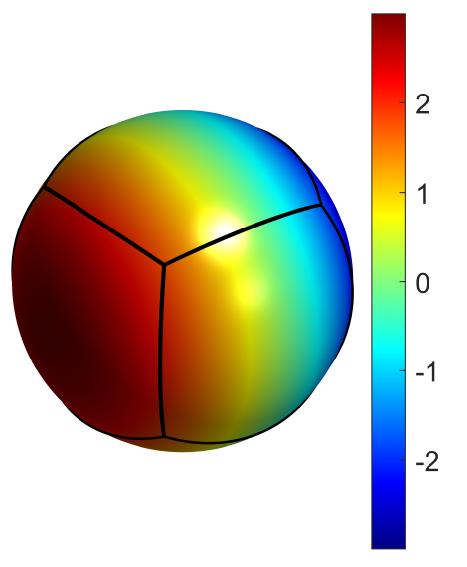}}\quad
\subfloat[]{\includegraphics[width = 0.22\textwidth, page = 8]{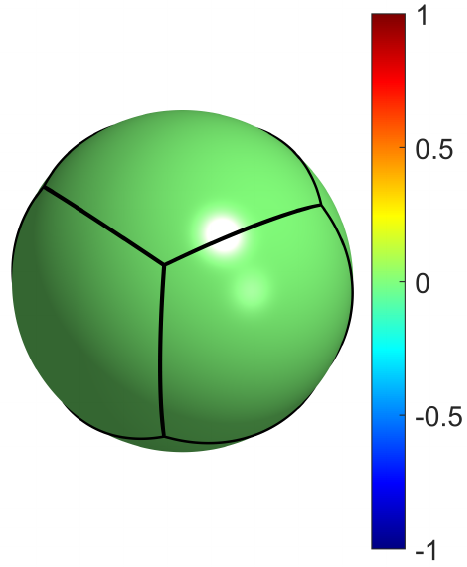}}\quad
\subfloat[]{\includegraphics[width = 0.22\textwidth, page = 9]{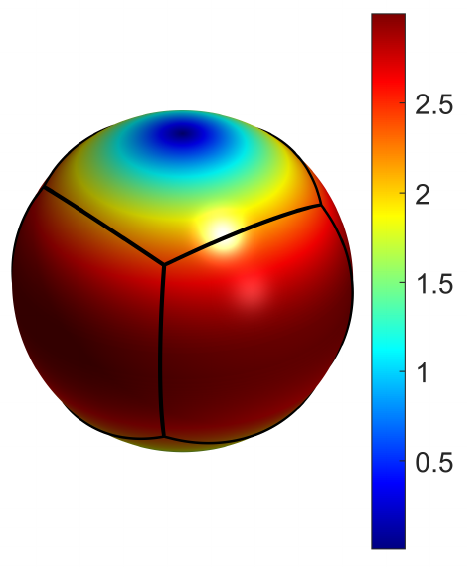}}
\caption{Components and magnitude of the exact surface traction field $\pmb{t}$ of the Stokes flow considered in Section \ref{sec:rotsphere}. In (a)--(c) the three components, in (d) the magnitude, i.e., $\norm{\pmb{t}}$. The colors and the colorbar identify the values pointwise. In particular the $z$ component is zero everywhere.}\label{fig:tractionrotsphere}
\end{figure}

\begin{figure}
\centering
\subfloat[]{\includegraphics[width = 0.22\textwidth, page = 10]{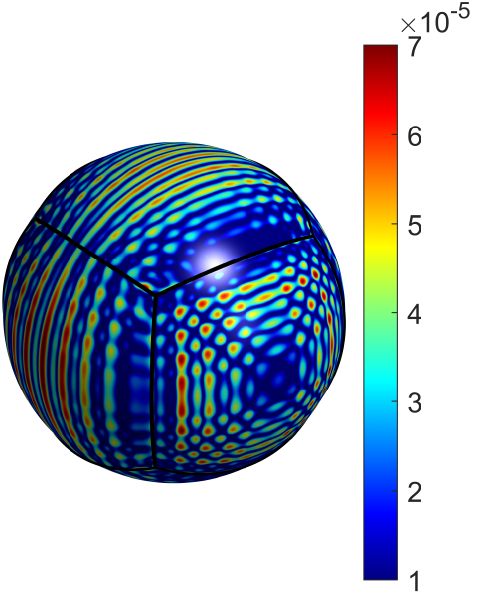}}\quad
\subfloat[]{\includegraphics[width = 0.22\textwidth, page = 11]{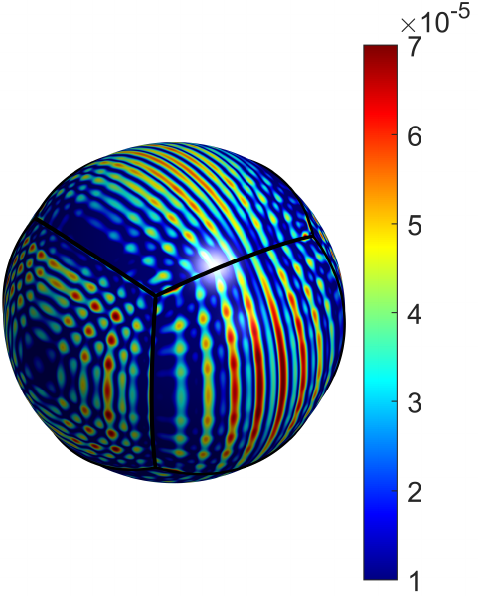}}\quad
\subfloat[]{\includegraphics[width = 0.22\textwidth, page = 12]{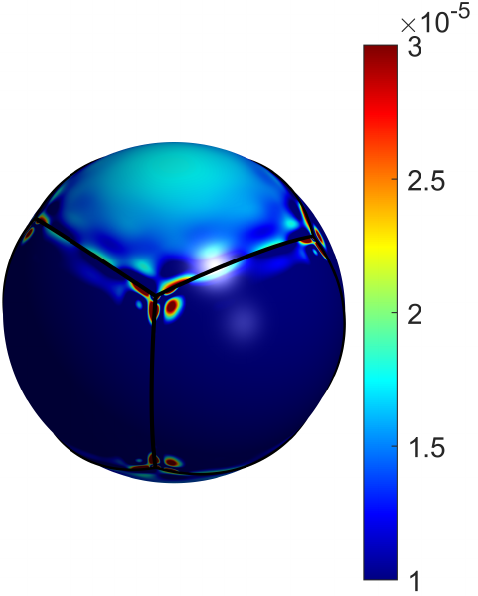}}\quad
\subfloat[]{\includegraphics[width = 0.22\textwidth, page = 13]{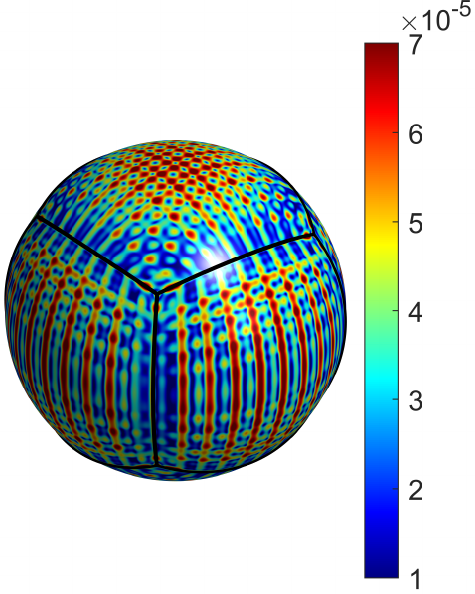}}
\caption{Pointwise error on the surface traction components and magnitude of the Stokes flow considered in Section \ref{sec:rotsphere}. The approximant spline space of each patch has degree $\m{d} = (2, 2)$ and it is defined on a mesh of $10 \times 10$ elements. In (a)--(c) the error on the three components, in (d) that on the magnitude of $\pmb{t}$. The colors and the color-bar identify the values pointwise.}\label{fig:tractionerrorrotsphere}
\end{figure}

In this setting, the velocity field and pressure of the flow are the solution of the Dirichlet Stokes problem
\begin{equation*}\label{eq:stokesdirichlet}
\left\{
\begin{array}{ll}\displaystyle -\eta \Delta \pmb{u}(\pmb{x}) + \nabla p(\pmb{x}) = \pmb{0} & \pmb{x} \in \RR^3\setminus S\\\\
\displaystyle\nabla \cdot \pmb{u}(\pmb{x}) = 0 & \pmb{x} \in \RR^3\setminus S\\\\
\pmb{u}|_{\Gamma}(\pmb{x}) = \omega \m{e}_3 \times \pmb{x} & \pmb{x} \in \Gamma
\end{array}
\right.
\end{equation*}
where $\m{e}_3 = (0, 0, 1)^T$. We can solve the problem by finding out the surface traction $\pmb{t}$ thorugh the BIE \eqref{eq:BIEmatrixform} and then using the representation formula \eqref{eq:birf} for computing velocity and pressure in the unbounded domain $\Omega = \RR^3\setminus S$. The expression of $\pmb{t}$ for this problem can be analytically computed \cite{chwang} and it is
$$
\pmb{t}(\pmb{x}) = - 3\eta\omega \frac{1}{R} \m{e}_3 \times \pmb{x}.
$$ 
Such expression shall be used here as reference for the computation of errors. Fixed $R = \omega = \eta = 1$ without loss of generality, the componentes of $\pmb{t}$ as well as its magnitude on the sphere $\Gamma$, written as 6 patch NURBS of rational degree 4, are reported in Figure \ref{fig:tractionrotsphere} (a)--(d). In particular we underline that $t_{\max} \defeq \max_{\pmb{x} \in \Gamma} \norm{\pmb{t}(\pmb{x})} = 3$.
We discretize the BIE using the IgA-BEM approach described in Section \ref{sec:igabem} and we compute the resulting integrals following the smoothly varying quadrature rule of Section \ref{sec:quadrule}. Here we consider spline spaces of degree $\m{d} = (2, 2)$ in all the patches for different mesh resolutions, namely meshes of $(2n) \times (2n)$ elements with $n = 2, \ldots, 5$. In Figure \ref{fig:tractionerrorrotsphere} we show the pointwise error for $n = 5$ on the traction components and magnitude. 

\begin{figure}
\centering
\subfloat[]{\includegraphics[width = 0.35\textwidth, page = 4]{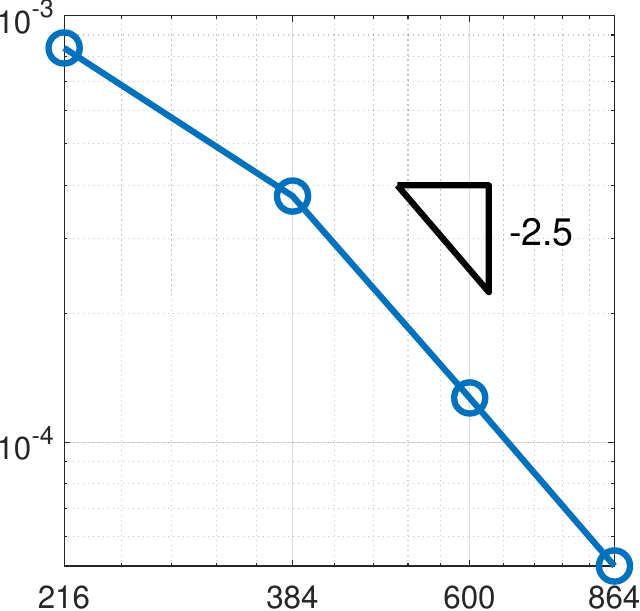}}\qquad
\subfloat[]{\includegraphics[width = 0.425\textwidth, page = 14]{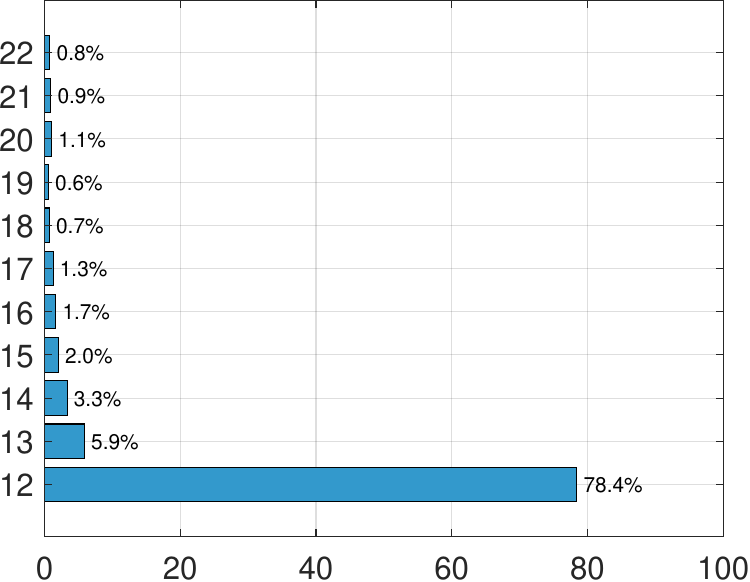}}
\caption{(a) $L^2$ error decay of the approximation of the surface traction and (b) number of nodes distribution for the numerical quadrature rule. In (a) the errore $e_t^{L^2}$ defined in \eqref{eq:meanL2error} is measured when we refine the mesh of $(2n) \times (2n)$ elements, for $n = 2, \ldots, 5$, corresponding to $216, 384, 600$ and $864$ degrees of freedom respectively. The reached rate of convergence is approximately $2.5$. In (b), for $n = 5$, we show the percentage of integrals approximated using a given number of quadrature nodes per direction. The $y$-axis reports the number of nodes, while the $x$-axis indicates the corresponding percentage. For readability, only numbers of nodes used in more than $0.5\%$ of the total integrals are shown on the $y$-axis. We highlight that $96.7\%$ of the integrals require fewer than 23 nodes, and the majority, about $78.4\%$, use the minimum number of nodes considered, namely 12. Much smaller percentages involve larger numbers of nodes.}\label{fig:L2errrotsphere}
\end{figure}

In order to quantify the overall quality of the discretizations and therefore study the decay of the approximation error, we follow the procedure of \cite{sauer} and define the relative traction error as
$$
e_t(\pmb{x}) \defeq \frac{\norm{\pmb{t}^h(\pmb{x}) - \pmb{t}(\pmb{x})}}{t_{\max}}
$$
and the $L^2$-norm of the error as
\begin{equation}\label{eq:meanL2error}
 e_t^{L^2} \defeq \sqrt{\frac{1}{|\Gamma|}\int_{\Gamma} e_t^2(\pmb{x})\dsigma}
\end{equation}
with $|\Gamma| = 4\pi R^2 = 4\pi$, the surface area of $\Gamma$. In Figure \ref{fig:L2errrotsphere} (a) we show the decay of $e_t^{L^2}$ against the number of degrees of freedom as we refine the mesh, i.e., for $n = 2, \ldots, 5$. As we can see, the achieved rate of convergence is approximately $2.5$, which outperforms the state-of-the-art \cite[Figure 23 (b)]{sauer} that stands out between $1$ and $1.5$. 

Finally in Figure \ref{fig:L2errrotsphere} (b), we analyze the distribution of the quadrature nodes. Namely, fixed $n = 5$, we show the percentage of integrals that use a given number of quadrature nodes, relative to the total number of integrals computed for the assembly of the system matrix $\bbol{S}$ of \eqref{eq:dirichlet}, that is, $864 \times 864 = 746.496$ integrals in total. In the $y$-axis of the hystogram we run over the numbers of quadrature nodes per direction, on the $x$-axis we have the corresponding percentages of integrals. Only number of nodes exceeding a $0.5\%$ percentage threshold are reported. For the chosen discretization settings, the minimal number of quadrature node per direction is $12$, which accounts for approximately $78.4\%$ of the entries in $\bbol{S}$. Numbers of nodes  greater than or equal to $23$ fall below the $0.5\%$ threshold and are thus omitted from the plot. Although the maximum number of nodes used is $41$, this occurs in fewer than $0.1\%$ of the integrals. This analysis is done in order to highlight that using a dynamic, distance-dependent, and support-specific number of nodes for the quadrature rule has a negligible impact on the overall computational cost, as only a very small percentage of integrals requires a significant increment of the nodes.

\subsection{Simulations of the exterior Stokes flow caused by a rising prolate spheroid}\label{sec:risingspheroid}
A prolate spheroid is a particular instance of ellipsoid in which the two axes in the equatorial plane are equal and the vertical axis is the major axis. Mathematically speaking, the coordinates of the surface points of a prolate spheroid verify the equation
$$
\frac{x_1^2 + x_2^2}{b^2} + \frac{x_3^2}{a^2} = 1
$$ 
for some $a, b \in \RR$ with $a > b$. The outward normal at point $\pmb{x} = (x_1, x_2, x_3)$ is therefore
\begin{equation}\label{eq:spheroidn}
\m{n}(\pmb{x}) = \nicefrac{\tilde{\m{n}}(\pmb{x})}{\norm{\tilde{\m{n}}(\pmb{x})}}\quad\text{with }\tilde{\m{n}}(x_1, x_2, x_3) \defeq \left(\frac{x_1}{b^2}, \frac{x_2}{b^2}, \frac{x_3}{a^2}\right).
\end{equation}
Defined $c \defeq \sqrt{a^2 - b^2}$, the surface area of a prolate spheroid is given by the closed formula \cite{spheroid}
\begin{equation}\label{eq:spheroida}
2\pi b^2 + 2\pi \frac{a^2b}{c} \arcsin\left(\frac{c}{a}\right).
\end{equation}
Assume now $S$ to be a solid prolate spheroid, with surface $\Gamma = \partial S$, immersed in a Newtonian fluid with dynamic viscosity $\eta$. The spheroid is rising in the fluid along the $x_3$-axis with velocity module $\omega$. This causes a creeping flow motion in the fluid.  In this settings, the velocity field and pressure of the flow are the solution of the following Dirichlet Stokes problem 
\begin{equation*}\label{eq:stokesdirichletspheroid}
\left\{
\begin{array}{ll}\displaystyle -\eta \Delta \pmb{u}(\pmb{x}) + \nabla p(\pmb{x}) = \pmb{0} & \pmb{x} \in \RR^3\setminus S\\\\
\displaystyle\nabla \cdot \pmb{u}(\pmb{x}) = 0 & \pmb{x} \in \RR^3\setminus S\\\\
\pmb{u}|_{\Gamma}(\pmb{x}) = -\omega \m{e}_3 & \pmb{x} \in \Gamma
\end{array}
\right.
\end{equation*}
where $\m{e}_3 = (0, 0, 1)^T$. We can solve the problem by finding out the surface traction $\pmb{t}$ through the BIE \eqref{eq:BIEmatrixform} and then using the representation formula \eqref{eq:birf} for computing velocity and pressure in the unbounded domain $\Omega = \RR^3\setminus S$. The expression of $\pmb{t}$ for this problem can be analytically computed \cite{spheroid1, spheroid2} and it is
\begin{equation}\label{eq:spheroidt}
\pmb{t}(\pmb{x}) = (4\pi a b^2)^{-1} \langle \m{n}(\pmb{x}), \pmb{x}\rangle \pmb{F}
\end{equation}
where $\m{n}$ is the outward normal defined in Equation \eqref{eq:spheroidn} and with $\pmb{F} = (0, 0, F_3)^T$ the \textbf{hydrodynamic drag} \cite[Section 4-30]{spheroid3}, i.e., a constant vector defined as follows. Set $\tau \defeq \nicefrac{a}{c}$ we have
$$
F_3 \defeq \frac{8\pi \eta c \omega}{(\tau^2 + 1) \arccoth(\tau) - \tau}
$$
where we recall that
$
\arccoth(\tau) = \frac{1}{2}\log\left(\frac{\tau + 1}{\tau - 1}\right).
$
The analytical expression of the traction \eqref{eq:spheroidt}  shall be used as reference for the computation of errors. In particular we underline that only the third component of $\pmb{t}$ is non-zero, because only the thrid component of the hydrodynamic drag $\pmb{F}$ is non-zero. 

\begin{figure}
\centering
\subfloat[]{\includegraphics[width = 0.22\textwidth, page = 16]{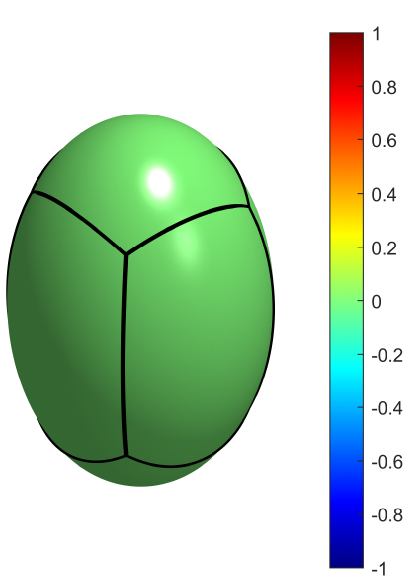}}\quad
\subfloat[]{\includegraphics[width = 0.22\textwidth, page = 17]{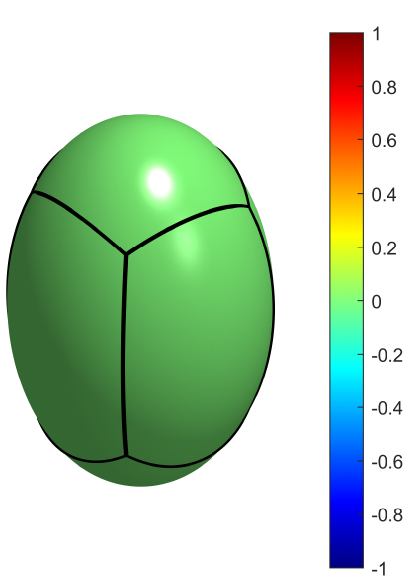}}\quad
\subfloat[]{\includegraphics[width = 0.22\textwidth, page = 18]{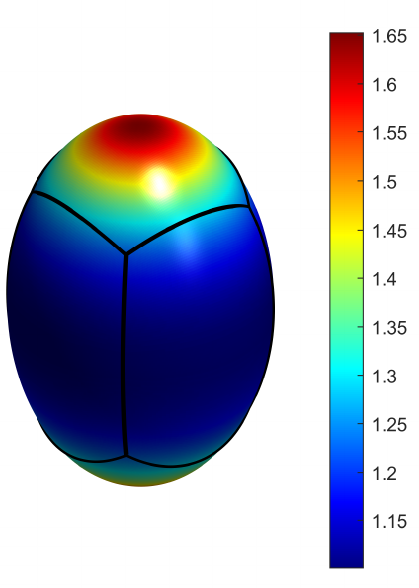}}\quad
\subfloat[]{\includegraphics[width = 0.22\textwidth, page = 19]{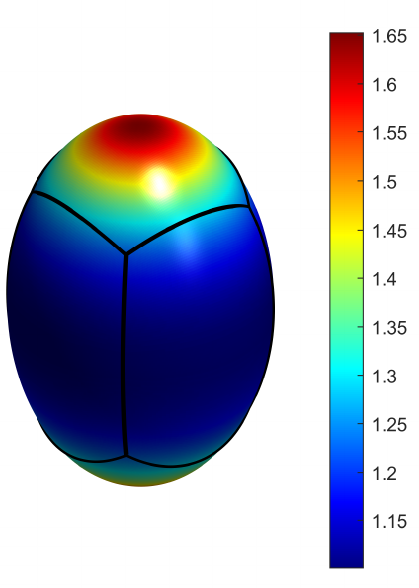}}
\caption{Components and magnitude of the exact surface traction field $\pmb{t}$ of the Stokes flow considered in Section \ref{sec:risingspheroid}. In (a)--(c) the three components, in (d) the magnitude, i.e., $\norm{\pmb{t}}$. The colors and the colorbars identify the values pointwise. We note in particular that $\norm{\pmb{t}} = t_3$ as the first two components of $\pmb{t}$ are zero and the third is positive everywhere. }\label{fig:tractionrotspheroid}
\end{figure}

\begin{figure}
\centering
\subfloat[]{\includegraphics[width = 0.22\textwidth, page = 20]{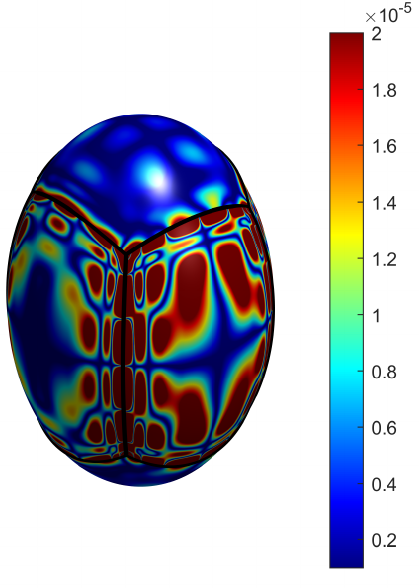}}\quad
\subfloat[]{\includegraphics[width = 0.22\textwidth, page = 21]{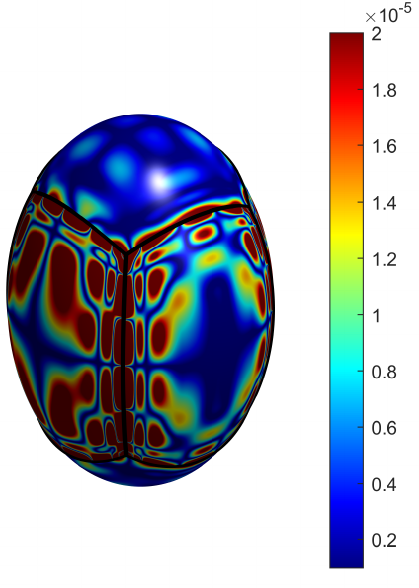}}\quad
\subfloat[]{\includegraphics[width = 0.22\textwidth, page = 22]{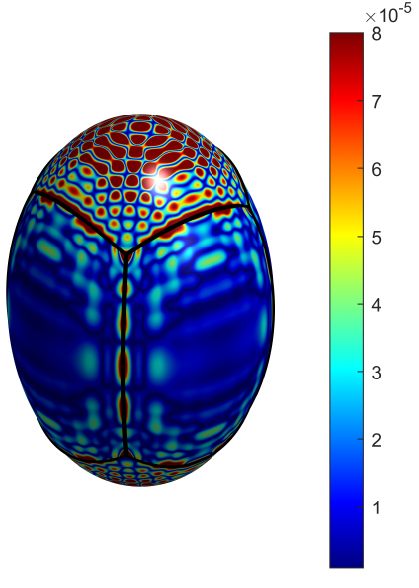}}\quad
\subfloat[]{\includegraphics[width = 0.22\textwidth, page = 23]{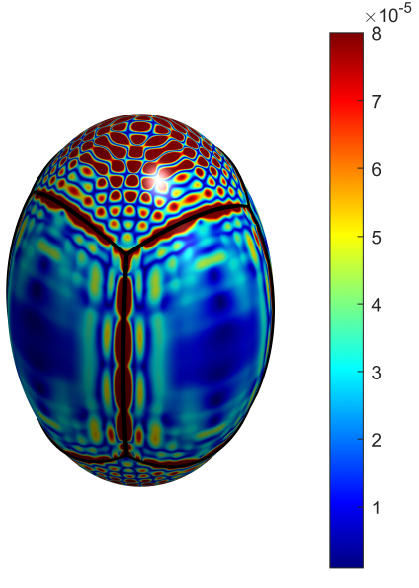}}
\caption{Relative pointwise error on the surface traction components and magnitude of the Stokes flow considered in Section \ref{sec:risingspheroid}. The approximant spline space of each patch has degree $\m{d} = (2, 2)$ and it is defined on a mesh of $10 \times 10$ elements. In (a)--(c) the error on the three components, in (d) that on the magnitude of $\pmb{t}$. The colors and the colorbars identify the values pointwise.}\label{fig:tractionerrorrotspheroid}
\end{figure}

\begin{figure}
\centering
\subfloat[]{\includegraphics[width = 0.35\textwidth, page = 15]{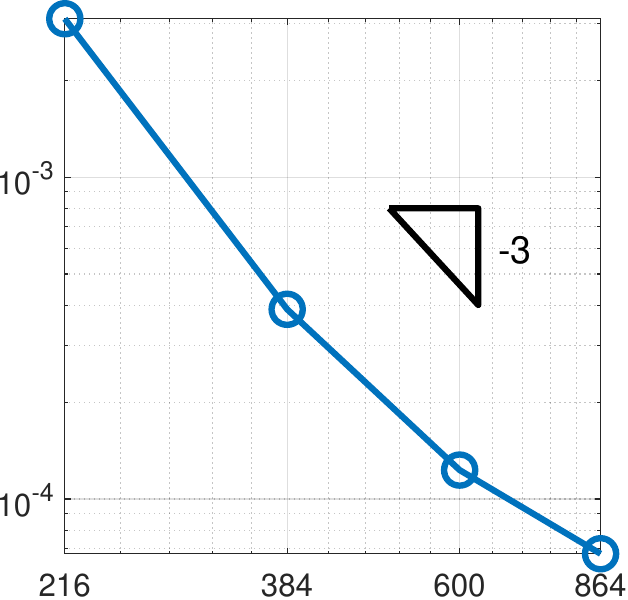}}\qquad
\subfloat[]{\includegraphics[width = 0.425\textwidth, page = 24]{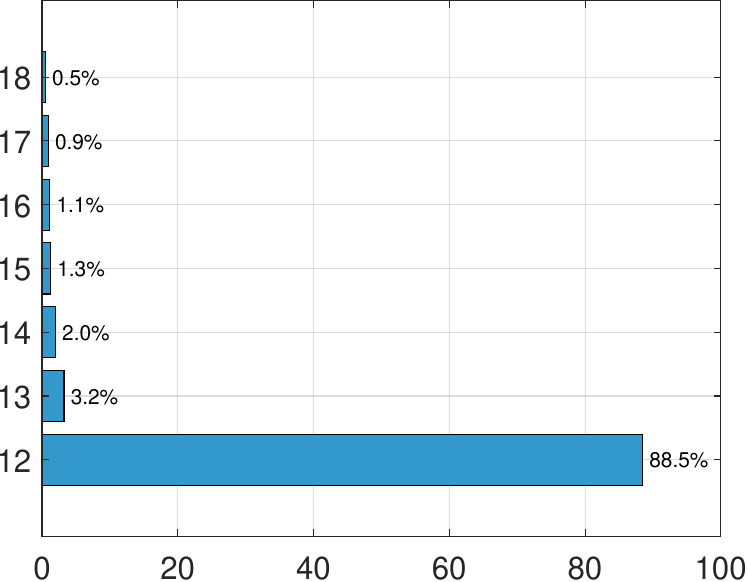}}
\caption{(a) $L^2$ error decay of the approximation of the surface traction and (b) number of nodes distribution for the numerical quadrature rule. In (a) the errore $e_t^{L^2}$ defined in \eqref{eq:meanL2error1} is measured when we refine the mesh of $(2n) \times (2n)$ elements, for $n = 2, \ldots, 5$, corresponding to $216, 384, 600$ and $864$ degrees of freedom respectively. The reached rate of convergence is approximately $3$. In (b), for $n = 5$, we show the percentage of integrals approximated using a given number of quadrature nodes per direction. The $y$-axis reports the number of nodes, while the $x$-axis indicates the corresponding percentage. For readability, only numbers of nodes used in more than $0.5\%$ of the total integrals are shown on the $y$-axis. We highlight that more than $97.5\%$ of the integrals use fewer than $19$ nodes, and the majority, about $88.5\%$, use the minimum number of nodes considered, namely $12$. Much smaller percentages involve larger numbers of nodes.}\label{fig:L2errrotspheroid}
\end{figure}

Fixed $a = 1.5$, $b = 1$ and $\omega = \eta = 1$, without loss of generality, the components of $\pmb{t}$ as well as its magnitude on the spheroid $\Gamma$, written as 6 patch NURBS of rational degree 4, are reported in Figure \ref{fig:tractionrotspheroid} (a)--(d). 
We discretize the BIE using the IgA-BEM approach described in Section \ref{sec:igabem} and we compute the resulting integrals following the smoothly varying quadrature rule of Section \ref{sec:quadrule}. We consider spline spaces of degree $\m{d} = (2, 2)$ in all the patches for different mesh resolutions, namely meshes of $(2n) \times (2n)$ elements with $n = 2, \ldots, 5$. In Figure \ref{fig:tractionerrorrotspheroid} we show the relative pointwise error for $n = 5$ on the traction components and magnitude. 
In order to quantify the overall quality of the discretizations and therefore study the decay of the approximation error, we define the relative traction error as
$$
e_t(\pmb{x}) \defeq \frac{\norm{\pmb{t}^h(\pmb{x}) - \pmb{t}(\pmb{x})}}{\norm{\pmb{t}(\pmb{x})}}
$$
and the $L^2$-norm of the error as
\begin{equation}\label{eq:meanL2error1}
 e_t^{L^2} \defeq \sqrt{\frac{1}{|\Gamma|}\int_{\Gamma} e_t^2(\pmb{x})\dsigma}
\end{equation}
with $|\Gamma|$ the surface area of the spheroid, that can be computed with the closed formula of Equation \eqref{eq:spheroida}. In Figure \ref{fig:L2errrotspheroid} (a) we show the decay of $e_t^{L^2}$ against the number of degrees of freedom as we refine the mesh, i.e., for $n = 2, \ldots, 5$. The achieved rate of convergence is approximately $3$, which outperforms the state-of-the-art \cite[Figure 27 (a)]{sauer} that stands out between $1$ and $1.5$. 

Finally in Figure \ref{fig:L2errrotspheroid} (b), we analyze the distribution of the quadrature nodes for this second test as well. Namely, fixed $n = 5$, we show the percentage of integrals that use a given number of quadrature nodes, relative to the total number of integrals computed for the assembly of the system matrix $\bbol{S}$ of \eqref{eq:dirichlet}, that is, $864 \times 864 = 746.496$ integrals in total. In the $y$-axis of the histogram we run over the numbers of quadrature nodes  in each direction, on the $x$-axis we have the corresponding percentages of integrals. Only number of nodes exceeding a $0.5\%$ percentage threshold are reported. For the chosen discretization settings, the minimal number of quadrature node per direction is $12$, which accounts for approximately $88.5\%$ of the entries in $\bbol{S}$. Numbers of nodes  greater than or equal to $19$ fall below the $0.5\%$ threshold and are thus omitted from the plot. The maximum number of nodes used is $27$. This occurs in about $0.02\%$ of the integrals.

\section{Conclusion}\label{sec:conclusion}
In this work, we have presented a new quadrature rule for isogeometric boundary element discretizations. The challenges posed by kernel singularities, namely, the integration of singular and nearly singular functions, as well as their distinction from regular ones, have been overcome through suitable changes of variables: the Duffy transformation for singular integrals and the Distance Calibrated Telles (DCT) transformation for non-singular ones. The latter is applied to all integrals whose domains do not include the polar point of the kernel, thus providing a unified treatment that removes the need to distinguish between nearly singular and regular integrals. This approach eliminates the sensitivity of the numerical accuracy to such classifications, which could otherwise lead to a significant loss of precision.
The proposed quadrature rule also relies on a support-wise direct spline integration formula, which reduces computational cost and memory access compared to classical element-wise integration, particularly for higher degree splines. The method has been validated through well established three-dimensional boundary element benchmarks for the Stokes problem in a multipatch setting, demonstrating superior performance compared to existing approaches.

As part of future work, we plan to extend the proposed procedure to spline spaces that enable local refinement. In the present formulation, the regularity between patches is $C^{-1}$, allowing non-conforming interfaces and some extent of non-uniformity. However, without local in-patch refinement, a non-negligible number of degrees of freedom may be wasted without improving accuracy. Incorporating local refinement would enable improved convergence both in cases where the solution exhibits localized features and in boundary element problems with regions of low regularity, such as edges and ridges. The formulation of the boundary integral equations in such situations is indeed more complex, and we expect that a dedicated treatment near these areas of reduced regularity will be required to preserve good convergence rates.

Finally, since the Stokes problem represents a fundamental step toward the simulation of the Navier–Stokes equations, we also plan to apply the proposed quadrature strategy within a boundary element framework for incompressible flow simulations in higher Raynold number regimes.

\section*{Acknowledgements}
The authors acknowledge the contribution of the National Recovery and Resilience Plan, Mission 4 Component 2 - Investment 1.4 - NATIONAL CENTER FOR HPC, BIG DATA AND QUANTUM COMPUTING (CUP B83C22002830001), funded by the European Union - Next Generation EU.\\
Francesco Patrizi aknowledges the MUR (Italian Ministry of University and Research) Excellence Department Project MatMod@TOV (CUP E83C23000330006) awarded to the Department of Mathematics of the University of Rome Tor Vergata. 

Cesare Bracco and Alessandra Sestini also acknowledge the partial support of MUR through the PRIN projects COSMIC (No. 2022A79M75) and NOTES (No. P2022NC97R), funded by the European Union - Next Generation EU.

The authors are members of Gruppo Nazionale per il Calcolo Scientifico - Istituto Nazionale di Alta Matematica (GNCS-INdAM). The INdAM support through GNCS 2025 projects \lq\lq{}PASTRAMI - sPline And Solver innovaTions foR Adaptive isogeoMetric analysIs\rq\rq{}, \lq\lq{}High-order BEM based numerical techniques for wave propagation problems\rq\rq{} and \lq\lq{}Tecnologie spline efficienti per l'analisi isogeometrica e l'approssimazione di dati\rq\rq{} (CUP E53C24001950001) is gratefully acknowledged. 

\bibliography{biblio}
\end{document}